\documentclass[10pt]{article}
\usepackage{color}
\usepackage{amsmath}
\usepackage{amssymb}
\usepackage{mathrsfs}
\usepackage[utf8]{inputenc}
\usepackage[english]{babel}
\usepackage{amsthm}
\usepackage{tikz}

\usepackage{graphicx}
\graphicspath{{./images/}}

\usepackage[toc,title,page]{appendix}

\usepackage{epsfig}

\usepackage[scaled=1]{helvet}
\usepackage{titlesec}
\usepackage{multido}
\usepackage{multirow}
\usepackage{blindtext}
\usepackage{color}
\usepackage{lipsum}
\usepackage{mathtools}
\usepackage[vcentermath]{youngtab}
\usepackage{young}
\usepackage[all,ps,dvips,graph]{xy}
\usepackage[misc,geometry]{ifsym}
\usepackage{moresize}

\usepackage[left=3cm,top=2cm,right=3cm,bottom=2cm]{geometry}
\usepackage{graphicx}
\numberwithin{equation}{subsection}

\let\OLDthebibliography\thebibliography
\renewcommand\thebibliography[1]{
  \OLDthebibliography{#1}
  \setlength{\parskip}{0pt}
  \setlength{\itemsep}{0pt plus 0.3ex}
}

\newtheorem{theorem}{Theorem}
\newtheorem{proposition}[theorem]{Proposition}
\newtheorem{corollary}[theorem]{Corollary}
\newtheorem{lemma}[theorem]{Lemma}

\newtheorem{example}[theorem]{Example}
\newtheorem{definition}[theorem]{Definition}
\newtheorem{conjecture}[theorem]{Conjecture}

\numberwithin{theorem}{section}

\newcommand{\bgamma}{\boldsymbol{\gamma}}

\newcommand{\bpartial}{\boldsymbol{\partial}}

\newcommand\bcalJ{\boldsymbol{\mathcal{J}}}

\newcommand\bcalNT{\boldsymbol{\mathcal{NT}}}
\newcommand\bcalNGAT{\boldsymbol{\mathcal{NT}}}
\newcommand\bcalD{\boldsymbol{\mathcal{D}}}

\newcommand\calA{\mathcal{A}}

\newcommand\calC{\mathcal{C}}
\newcommand\calD{\mathcal{D}}
\newcommand\calE{\mathcal{E}}
\newcommand\calF{\mathcal{F}}

\newcommand\calL{\mathcal{L}}

\newcommand\calNT{\mathcal{NT}}
\newcommand\ncalNT{\mathcal{NT}_n}

\newcommand\calP{\mathcal{P}}
\newcommand\calQ{\mathcal{Q}}

\newcommand{\F}{ { \mathbb F}}

\newcommand{\mSLTM}{{$m-$\rm SLTM}}

\newcommand{\nSLTM}{{$n-$\rm SLTM}}
\newcommand{\nSLTMn}{$\mathcal{TM}_n$}

\newcommand{\SLTM}{{\rm SLTM}}
\newcommand{\SLTMn}{$\mathcal{TM}$}

\newcommand{\ulu}{\underline{u}}
\newcommand{\ulv}{\underline{v}}
\newcommand{\ulw}{\underline{w}}

\newcommand{\ulz}{\underline{z}}

\begin{document}
  \Yvcentermath1
\title{The category $\bcalNT:$ Isomorphism criteria and applications.}
\author{\sc  Diego Lobos}

\maketitle


\begin{abstract}
  The category $\bcalNT$ was introduced in \cite{Lobos2} in order to provide a structural setting to the study of certain commutative algebras appearing in the context of the diagrammatic Hecke category of Elias and Williamson \cite{EW}. The category $\bcalNT$ have as objects, all that we called \emph{nil graded algebras associated to triangular matrices} and as morphisms, all the \emph{preserving degree} homomorphisms of graded algebras between them. In this article we develop a series of \emph{isomorphism criteria} on $\bcalNT$ and we show some applications relative to the search of a general classification and the particular case of Gelfand-Tsetlin subalgebras coming from Soergel Calculus.
\end{abstract}
keywords: Commutative graded algebras, Jucys-Murphy elements, Soergel Calculus.
\tableofcontents

\pagenumbering{arabic}

\section{Introduction}
\subsection{Motivation}
The \emph{diagrammatic Hecke Category} $\bcalD,$ is a linear and monoidal category, it was defined in 2016 by Elias and Williamson \cite{EW} in order to provide a diagrammatic approach to the study of the  category of \emph{Soergel bimodules,} a well-known incarnation of the \emph{Hecke category,} \cite{So}. Since then, the category $\bcalD$ have gained attention of both algebraists and algebraic geometrist, due to its intrinsic connection with both branches of study. For further studies on this topic, we recommend to see \cite{EW-intro} and \cite{LibGentle} for example.

A particular interest is the study of the internal structure of the category $\bcalD$ itself, that is, the structure of the many $\textrm{\bf Hom}$ spaces and more particularly of the many $\textrm{\bf End}$ algebras. In this context, for example, the diagrammatic setting provided by Elias and Williamson have been crucial to find out and prove the cellularity of the \emph{Light leaves bases} of Libedinsky \cite{LibLightLeaves},\cite{EW}. Another example of the power of this diagrammatic approach is the discovery due to Ryom-Hansen, of explicit \emph{Jucys-Murphy} elements relative to the Libedinsky's bases \cite{Steen2} (In the following when we mention Jucys-Murphy elements in $\bcalD,$ we refer to those developed by Ryom-Hansen in \cite{Steen2}). Finally we would like to mention, probably the most important application of this diagrammatic setting, the disproof of the famous Lusztig's conjecture, given by Williamson in \cite{Williamson2017}.

In \cite{LiPl}, Libedinsky and Plaza, found out a connection between the category $\bcalD$ and what they called the \emph{Blob category}, another diagrammatic category, inspired in the well-known \emph{KLR-}algebras of Khovanov, Lauda and Rouquier. The similar patterns that they observed in both categories motivate them to define the \emph{Blob vs Soergel conjecture}, where roughly speaking, they declare that these two categories are equivalent. Some of the first  evidences on the veracity of the aforementioned conjecture can be found in the work of Espinoza and Plaza \cite{Esp-Pl}, where explicit isomorphisms were built between the subalgebras generated by the Jucys-Murphy elements, coming on one side from \emph{Two-color} Soergel Calculus and on the other of certain idempotent truncation of the KLR version of the \emph{blob algebra} (see \cite{PlazaRyom}). In \cite{Lobos-Plaza-Ryom-Hansen}, Plaza, Ryom-Hansen and the author, obtain isomorphisms in a higher level, this time between endomorphism algebras of $\calD$ in type $\widetilde{A}_1$ and the corresponding idempotent truncation of blob algebras. In both cases the procedure implies the search of adequate presentations in terms of generators and relations for the algebras involved, therefore, in both \cite{Esp-Pl} and \cite{Lobos-Plaza-Ryom-Hansen}, the reader will find, not only a clear proof of each isomorphism, but also more information on the structure of the corresponding algebras.

Some years ago, Bowman, Cox and Hazi have proved in \cite{bow-cox-hazi}, the \emph{Blob vs Soergel} conjecture, they in fact, proved a more general result, obtaining as a consequence the equivalences predicted by Libedinsky and Plaza.  The procedure used in \cite{bow-cox-hazi}, based in alcove geometry and paths configurations,  allows the authors to cover all the cases in just one go, while in \cite{Esp-Pl} and  \cite{Lobos-Plaza-Ryom-Hansen} only particular cases were covered. On the other hand the general approach developed in \cite{bow-cox-hazi} do not implies the aforementioned internal information on the corresponding algebras obtained in \cite{Esp-Pl} and \cite{Lobos-Plaza-Ryom-Hansen} or any generalization of that. Therefore, we consider as an open problem the search of this internal information for all the cases no covered in \cite{Esp-Pl} and \cite{Lobos-Plaza-Ryom-Hansen}, that is, in principle, to find out presentations for all the algebras related with the \emph{Blob vs Soergel} conjecture.

In the last few years, we have been developing certain strategies to address this problem. For example in \cite{Lobos1}, we obtained presentations for the subalgebras generated by Jucys-Murphy elements in certain family of idempotent truncations of the \emph{generalized blob algebras}. Later in \cite{Lobos2}, working on the side of the category $\bcalD,$ we developed a series of techniques and procedures to obtain  presentations for the $JM-$\emph{algebras}, that is, the subalgebras generated by Jucys-Murphy elements (also called Gelfand-Tsetlin subalgebras), in the many $\textrm{\bf End}$ algebras for any type (far beyond the initial objective: type $\widetilde{A}$). Both \cite{Lobos1} and \cite{Lobos2} are natural generalizations of \cite{Esp-Pl}.

In general the subalgebra generated by the \emph{Jucys-Murphy} elements in a cellular algebra $\mathcal{A},$ is a commutative subalgebra that under certain conditions, it could be a maximal commutative subalgebra of $\calA$  (see \cite{Mat-So}). One can expect that $JM-$subalgebra contains important information on the algebra $\mathcal{A}$ itself. On the other hand \emph{Jucys-Murphy} elements are very relevant in the study of the \emph{Representation Theory} of the corresponding cellular algebra $\calA$ (see \cite{Mat-So}). We think that all the study that we are doing on $JM-$algebras in the context of the category $\bcalD,$ will help to learn more on this important category.

The category $\bcalNT$ was introduced in \cite{Lobos2} in order to provide a structural setting to the study of the many $JM-$algebras appearing in the context of the diagrammatic Soergel category $\bcalD.$ The category $\bcalNT$ have as objects, all that we called \emph{nil graded algebras associated to triangular matrices} and as morphisms, all the \emph{preserving degree} homomorphisms of graded algebras between them. Naturally, the aforementioned $JM-$algebras of $\bcalD$ are (the main) objects in our category, and all the knowledge of $\bcalNT$ developed in \cite{Lobos2} allowed us to find out the desired presentations for those algebras. In particular we would like to use our work in the search of a complete classification of the many $JM-$algebras of the Soergel category $\bcalD.$

On the other hand, the category $\bcalNT$ is, in our opinion, by itself an interesting object of study, since their objects are basically certain quotients of polynomial algebras on several variables, they have potential connections with other branches of study, for example with Algebraic Geometry. To obtain tools to classify the objects on $\bcalNT,$ could help, for example, to obtain new tools for analysis of certain types of Algebraic Varieties.

Being more specific, the objects of the category $\bcalNT$ are commutative graded algebras defined by a presentation with a finite set of generators and certain quadratic relations codified by a strictly lower triangular matrix. The dependence on the \emph{codifying} matrices, invite us to ask how those matrices, could be relevant when we want to know if two of those algebras are isomorphic or not. Moreover, we would like to know if a complete classification of the objects of $\bcalNT$ is possible, and which role play the corresponding codifying matrices in these aspects. In this article we give some first steps in the search of answers for those questions.  The content of this article is organized as follows:
\begin{itemize}
  \item We start Section \ref{sec-Isom-criteria} by recalling some definitions and first results on the category $\bcalNT$ coming from \cite{Lobos2}.
  \item In subsection \ref{ssec-Fund-eq}, we obtain our first isomorphism criterion. In Theorem \ref{coro-Key-condition} we provide an absolute (but not easy to apply) matrix criterion to determine if two algebras in $\bcalNT$ are or not isomorphic. Theorem \ref{coro-Key-condition}, basically asserts that two algebras $\calA(T)$ and $\calA(S)$ in the category $\bcalNT,$ are isomorphic if and only if there exist an invertible matrix $\Gamma,$ whose entries satisfy certain  system of quadratic equations, determined by the entries of $T$ and $S.$
  \item In subsection \ref{ssec-ETO}, we define a set of matrix operations, that we called \emph{Elementary triangular operations} (or ETO's for shortness). This matrix operations act on the set of triangular matrices, each one subject to certain restrictions. In Definition \ref{def-ETO-assignments}, we define \emph{Equivalence by ETO's}, an equivalence relation between triangular matrices. In Theorem \ref{theo-sim-vs-approx}, we show that if two matrices $T,S$ are equivalent by ETO's, then the corresponding algebras $\calA(T),\calA(S)$ are isomorphic as objects of the category $\bcalNT.$
  \item In Proposition \ref{prop1-sim-vs-approx}, we connect the matrix $\Gamma$ coming from Theorem \ref{coro-Key-condition} and ETO's. We prove that if $\Gamma$ is any elementary matrix, defining an isomorphism $\gamma:\calA(T)\rightarrow \calA(S),$ then it implies that the matrices $T$ and $S$ are equivalent by ETO's. Since any isomorphism $\calA(T)\rightarrow \calA(S)$ is determined by an invertible matrix (Theorem \ref{coro-Key-condition}), and any invertible matrix $\Gamma$ can be factored as a product of elementary matrices, it seem sensible to think that the existence of an isomorphism between $\calA(T)$ and $\calA(S),$ will imply that the matrices are necessarily equivalent by ETO's. We declare this in Conjecture \ref{conj-sim-vs-approx}. We do not have a proof of this yet, but we have an important set of evidences on its veracity, we show some of them in this article.
      \item In Section \ref{sec-int-applic} we show some applications of the criteria obtained in Section \ref{sec-Isom-criteria}. In Subsection \ref{ssec-small-cases}, we obtain a complete classification of algebras $\calA(T),$ where $T$ is a lower matrix of size $2$ or $3$ (see Propositions \ref{Prop-Nn-for-n2} and \ref{prop-Cl-eq-2-for-n3}). We also prove for those cases the veracity of Conjecture \ref{conj-sim-vs-approx} (see Corollaries \ref{coro-Prop-Nn-for-n2} and \ref{coro-Prop-Nn-for-n3}).
      \item In Subsection \ref{ssec-class-of-0n}, we use Theorem \ref{coro-Key-condition} to deduce an absolute criterion, depending explicitly on the entries of a matrix $U,$ to determine if the algebras $\calA(U)$ and $\calA(0_n)$ are isomorphic or not. Theorem \ref{theo-zero-class-n-gen}, is a clear example of what we were talking about, when we propose ourselves to look for criteria coming from the codifying matrices, to determine when two algebras in $\bcalNT$ are isomorphic or not. Although it only cover one isomorphism class, we consider this result as a challenging invitation to obtain such as explicit criteria for other isomorphism classes. This is the motivation to keep on working in the search of a complete classification of the objects in $\bcalNT.$ On the other hand Proposition \ref{prop2-sim-vs-approx} provide another evidence on the veracity of Conjecture \ref{conj-sim-vs-approx}. In this case, Proposition \ref{prop2-sim-vs-approx} asserts that if an algebra $\calA(U)$ is isomorphic to $\calA(0_n),$ then necessarily the matrices $U$ and $0_n$ are equivalent by ETO's.
      \item As we mentioned above, in Subsection \ref{ssec-small-cases}, we obtained a complete classification for algebras codified by matrices of size $n=2$ or $n=3.$ In fact, Propositions \ref{Prop-Nn-for-n2} and \ref{prop-Cl-eq-2-for-n3}, provide explicitly the number of classes $N_n$ for those sizes. In both cases, the number is finite and it is independent on the ground field. Those facts invite us to open the following questions:
          \begin{description}
            \item[Q1:]  Is the number of classes always finite?
            \item[Q2:] Is the number $N_n$ always independent of the ground field?
          \end{description}
      \item Independently on if the number $N_n$ is finite or not, in Subsection \ref{ssec-lower-bound}, we obtain a lower bound for the number of classes for each size $n,$ that is independent on the ground field. In Theorem \ref{theo-lower-bound}, we describe an explicit set of representatives of different classes. This result will be useful, when we have answered \textbf{Q1} and \textbf{Q2}, and naturally look for the exact value of $N_n.$
      \item At the end of this article, we show some first applications to the problem of classifying $JM-$algebras appearing in the context of the diagrammatic Soergel category $\bcalD.$
\end{itemize}

The main results of this article are:
\begin{itemize}
  \item Theorem \ref{coro-Key-condition}, where we obtain a system of quadratic equations, that determine in terms of the codifying matrices, if two algebras in $\bcalNT$ are isomorphic or not.
  \item Theorem \ref{theo-sim-vs-approx}, where we provide a way to navigate on an isomorphism class, via our elementary triangular operations.
  \item Theorem \ref{theo-zero-class-n-gen}, where, as an application of Theorem \ref{coro-Key-condition}, we describe explicitly all the element of a particular isomorphism class in $\bcalNT.$
  \item Theorem \ref{theo-lower-bound}, where as application of Theorem \ref{coro-Key-condition}, we obtain a first lower bound for the number of classes for each size $n.$
  \item Theorem \ref{theo-appl-commuting-moves}, where we prove that the $JM-$algebra $\calA_{\ulu},$ is invariant under commuting moves and Proposition \ref{prop-calA-not-invariant-under-braid-moves}, where we show that the algebra $\calA_{\ulu}$ is not invariant under other braid moves.
\end{itemize}

\subsection{Notations and terminology}\label{ssec-notations}
The following notation will be used along this article:
\begin{itemize}
  \item $\F$ is a fixed field of characteristic different from $2.$ We also denote $\F^{\times}=\F-\{0\}.$
  \item $0_n$ is the zero matrix and $I_n$ is the identity matrix, both of size $n\times n$
  \item We denote by $B_{n,l}=[b_{ij}], (1\leq l< n),$ the  $n\times n-$matrix given by
  \begin{equation}\label{def-Snl-matrix}
    b_{ij}=\left\{\begin{array}{cc}
                    1 & \textrm{if}\quad i=n\quad\textrm{and}\quad 1\leq j\leq l \\
                    \quad & \quad \\
                    0 & \textrm{otherwise}
                  \end{array}\right.
  \end{equation}
  \item A \nSLTM,  is a strictly lower $n\times n-$matrix,  with entries in $\F.$ A \SLTM, is a strictly lower matrix (of any size), with entries in $\F.$ We denote by \nSLTMn, the set of all the \nSLTM,  and by \SLTMn,  the set of all the \SLTM.
  \item If $T$ is a \nSLTM, and $1\leq k<n,$ then we denote by $T|_k$ the matrix obtained from $T$ by erasing the rows and the columns of $T$ from the $(k+1)-$th to the $n-$th.
  \item Given a \nSLTM, $U=[u_{ij}],$ an scalar $\alpha\in\F$ and indices $1\leq i<j<k\leq n,$ we denote by $\Delta^{(\alpha)}_{i,j,k}(U)$ the element $\Delta^{(\alpha)}_{i,j,k}(U):=\alpha u_{ki}+u_{kj}u_{ji}\in \F.$
  When there is no possible confusion we only write $\Delta^{(\alpha)}_{i,j,k}$ instead of $\Delta^{(\alpha)}_{i,j,k}(U).$
  \item Given a \nSLTM, $T,$  a \mSLTM, $S$ and a $m\times n-$matrix $C,$ we define the matrix $T\nabla_{C}S$ by blocks, as follows:
    \begin{equation*}
      T\nabla_{C}S=\left[\begin{matrix}
                           T & 0 \\
                           C & S
                         \end{matrix}\right]
    \end{equation*}
\end{itemize}

\section{Isomorphism criteria}\label{sec-Isom-criteria}

\subsection{Generalities}
In this subsection we recall the main facts on the category $\bcalNT$ introduced in \cite{Lobos2}.

\begin{definition}\label{Def-A(T)}
  Given a \nSLTM, $T=[t_{ij}],$ the \emph{nil graded algebra associated to} $T$ is the commutative unitary $\F-$algebra defined by generators $X_1,\dots,X_n$ and relations:
  \begin{equation}\label{EQ-def-A(T)}
   { \left\{\begin{array}{c}
      X_1^{2}=0 \\
       \quad\\
      X_i^{2}={\sum_{j<i}}{t_{ij}}X_jX_i,\quad (2\leq i \leq n)
    \end{array}\right.}
  \end{equation}
  { We} denote this algebra by $\calA(T).$  We define a grading on $\calA(T),$ by $\deg(X_i)=2,$ for any $1\leq i\leq n.$
\end{definition}

Note that Definition \ref{Def-A(T)} implies that the algebra $\calA(T)$ is a \emph{nil graded algebra}, that is, a graded algebra, where any homogeneous element with degree different from $0,$ is a nilpotent element (see \cite{Lobos2} for details). In particular, if $T=0_n$ then Definition \ref{Def-A(T)} implies that $\calA(T)$ is generated by square nilpotent elements.

We denote by $\ncalNT$ the family of all nil graded algebras defied by some \nSLTM, and $\calNT$ the family of all nil graded algebras defined by some \SLTM.

The following Theorem was proved in \cite{Lobos2}:

\begin{theorem}\label{theo-monomial-basis}
  If $T=[t_{ij}]$ is a \nSLTM, then the dimension of the algebra $\calA(T)$ is $2^n$ and the set $\{X_1^{a_1}\cdots X_n^{a_n}:a_i\in\{0,1\}\}$ is a basis for $\calA(T).$
\end{theorem}

\begin{definition}\label{def-preserving-degree-homo}
  Given two \SLTM, $T,S,$ a \emph{preserving degree homomorphism}  $\gamma:\calA(T)\rightarrow\calA(S),$ is any homomorphism of algebras such that $\deg(X_i)=\deg(\gamma(X_i))=2,$ for any $1\leq i\leq n.$ If also $\gamma$ is an isomorphism we say that $\gamma$ is a \emph{preserving degree isomorphism}.
\end{definition}


\begin{definition}\label{def-category-bcalNT}
  We define the category $\bcalNT$  as the one whose set of objects is $\calNT$ and whose morphisms are all the \emph{preserving degree homomorphisms} between them. An isomorphism in $\bcalNT$ is a \emph{preserving degree isomorphism} between two members of $\calNT.$
\end{definition}

In the following we shall say that the algebras $\calA(T)$ and $\calA(S)$ are isomorphic if they are isomorphic as object of the category $\bcalNT$ (that is, in the sense of Definition  \ref{def-category-bcalNT}).

In this article we are interested in to find out certain matrix criteria, to determine when two objects $\calA(T),\calA(S)$ of $\bcalNT,$ are isomorphic. Given $T,S$ two \nSLTM,  we denote by $T\sim S$ whenever the algebras $\calA(T),\calA(S)$ are isomorphic. Relation $\sim$ is a equivalence relation on the set of all \nSLTM.


\subsection{Fundamental equation}\label{ssec-Fund-eq}
We shall assume that $T=[t_{ij}]$ is a \nSLTM,  $S=[s_{ij}]$ is a $m-$\SLTM, and $X_j,Y_i$ are the generators defining $\calA(T)$ and $\calA(S)$ respectively.

Note that any morphism $\gamma:\calA(T)\rightarrow\calA(S),$ defines a matrix $\Gamma=[\gamma_{ij}],$ by the equation
\begin{equation}\label{EQ-def-matrix-Gamma}
  \gamma(X_j)=\sum_{i=1}^{m}{\gamma_{ij}}Y_i,\quad(1\leq j\leq n)
\end{equation}
In particular any morphism is totally determined by a \emph{linear transformation} $\textrm{spam}\{X_1,\dots,X_n\}\rightarrow \textrm{spam}\{Y_1,\dots,Y_m\}.$ On the other hand, a given linear transformation $\textrm{spam}\{X_1,\dots,X_n\}\rightarrow \textrm{spam}\{Y_1,\dots,Y_m\},$ not necessarily defines a morphism in the category $\bcalNT.$ The following lemma provides a first matrix criterion to find out if a given matrix $\Gamma$ defines a morphism $\gamma:\calA(T)\rightarrow\calA(S)$ via equation \ref{EQ-def-matrix-Gamma}, or not:

\begin{lemma}\label{lemma-Key-condition}
  Given a $m\times n-$ matrix, $\Gamma=[\gamma_{ij}],$ then the assignment
  \begin{equation*}
    X_j\mapsto \sum_{i=1}^{m}{\gamma_{ij}}{Y_i}
  \end{equation*}
  defines a preserving degree homomorphisms of algebras $\gamma:\calA(T)\rightarrow\calA(S),$ if and only if
  \begin{equation}\label{Key-EQ}
    2\gamma_{ir}\gamma_{kr}+\gamma_{kr}^{2}s_{ki}=
    \sum_{j<r}{t_{rj}}(\gamma_{kj}\gamma_{kr}s_{ki}+\gamma_{kj}\gamma_{ir}+\gamma_{ij}\gamma_{kr}),\quad (1\leq r\leq n;1\leq i<k\leq m).
  \end{equation}
\end{lemma}

\begin{proof}
   The assignment $X_r\mapsto\bgamma(X_r)= \sum_{i=1}^{m}\gamma_{jr}Y_j,$ defines a morphism $\bgamma:\calA(T)\rightarrow\calA(S)$ in the category $\bcalNGAT$ if and only if for each $r=1,\dots,n$ we have:
  \begin{equation}\label{eq1-proof-lemma-conditions-for-morphisms}
    (\gamma(X_r))^2=\sum_{j<r}{t_{rj}\gamma(X_j)\gamma(X_r)}.
  \end{equation}

  If we develop the left side of equation \ref{eq1-proof-lemma-conditions-for-morphisms}, we obtain:
  \begin{equation*}
    (\gamma(X_r))^2=\left(\sum_{i=1}^{m}\gamma_{ir}Y_i\right)^2=\sum_{k=2}^{m}\sum_{i<k}{2\gamma_{ir}\gamma_{kr}Y_iY_k}
    +\sum_{k=1}^{m}\gamma_{kr}^2Y_k^2.
  \end{equation*}

 Now by relation \ref{EQ-def-A(T)}, we obtain
 \begin{equation*}
    (\gamma(X_r))^2=
    \sum_{k=2}^{m}\sum_{i<k}\left(2\gamma_{ri}\gamma_{kr}+\gamma_{kr}^2s_{ki}\right){Y_iY_k}
  \end{equation*}

On the other hand, if we develop the expression on the right in equation \ref{eq1-proof-lemma-conditions-for-morphisms}, we obtain:

\begin{equation}\label{eq3-proof-lemma-conditions-for-morphisms}
  \sum_{j<r}{t_{rj}\gamma(X_j)\gamma(X_r)}=
  \sum_{j<r}{t_{rj}\left(\sum_{k=1}^{m}\gamma_{kj}Y_k\right)\left(\sum_{i=1}^{m}\gamma_{ir}Y_i\right)}
\end{equation}
If we develop the right side of equation \ref{eq3-proof-lemma-conditions-for-morphisms}, we obtain

\begin{equation*}
  \sum_{j<r}{t_{rj}\gamma(X_j)\gamma(X_r)}=
  \sum_{j<r}{t_{rj}\left(\sum_{k=2}^{m}\sum_{i<k}\left(\gamma_{kj}\gamma_{ir}
  +\gamma_{ij}\gamma_{kr}\right)Y_iY_k+
  \sum_{k=1}^{m}\gamma_{kj}\gamma_{kr}Y_k^2\right)}
\end{equation*}

Again by relation \ref{EQ-def-A(T)}, we obtain

\begin{equation*}
  \sum_{j<r}{t_{rj}\gamma(X_j)\gamma(X_r)}=
  \sum_{j<r}{t_{rj}\sum_{k=2}^{m}\sum_{i<k}\left(\gamma_{kj}\gamma_{ir}+\gamma_{ij}\gamma_{kr}+
  s_{ki}\gamma_{kj}\gamma_{kr}\right)Y_iY_k}
\end{equation*}

or equivalently:

\begin{equation}\label{eq5-2-proof-lemma-conditions-for-morphisms}
  \sum_{j<r}{t_{rj}\gamma(X_j)\gamma(X_r)}=
 {\sum_{k=2}^{m}\sum_{i<k} \sum_{j<r}t_{rj}\left(\gamma_{kj}\gamma_{ir}+\gamma_{ij}\gamma_{kr}+
  s_{ki}\gamma_{kj}\gamma_{kr}\right)Y_iY_k}
\end{equation}

Equations \ref{eq1-proof-lemma-conditions-for-morphisms}, \ref{eq5-2-proof-lemma-conditions-for-morphisms} and the linear independence of the monomial basis $\{Y_1^{a_1}\cdots Y_m^{a_m}:a_i\in\{0,1\}\}$ (Theorem \ref{theo-monomial-basis}), imply that:

\begin{equation*}
   2\gamma_{ir}\gamma_{kr}+\gamma_{kr}^2s_{ki}=\sum_{j=1}^{n} t_{rj}\left(\gamma_{kj}\gamma_{kr}s_{ki}+\gamma_{kj}\gamma_{ir}+\gamma_{ij}\gamma_{kr}\right).
  \end{equation*}
  for any $r\in\{1,\dots,n\}$ and any pair $i,k\in\{1,\dots,m\}$ such that $i<k,$ as desired.
\end{proof}

\begin{theorem}\label{coro-Key-condition}
Let $T,S$ be two \nSLTM. Given a $n\times n-$ matrix, $\Gamma=[\gamma_{ij}],$ then the assignment
  \begin{equation*}
    X_j\mapsto \sum_{i=1}^{n}{\gamma_{ij}}{Y_i}
  \end{equation*}
  defines an isomorphisms $\bgamma:\calA(T)\rightarrow\calA(S),$ if and only if  $\Gamma$ is an invertible matrix whose entries satisfy the following system of equations:
  \begin{equation}\label{Key-EQ-coro}
    2\gamma_{ir}\gamma_{kr}+\gamma_{kr}^{2}s_{ki}=
    \sum_{j<r}{t_{rj}}(\gamma_{kj}\gamma_{kr}s_{ki}+\gamma_{kj}\gamma_{ir}+\gamma_{ij}\gamma_{kr}),\quad (1\leq r\leq n;1\leq i<k\leq n).
  \end{equation}
\end{theorem}

\begin{proof}
If $\bgamma$ is an isomorphism, then clearly the associated linear transformation (obtained by restriction) $\textrm{spam}\{X_1,\dots,X_n\}\rightarrow \textrm{spam}\{Y_1,\dots,Y_m\},$ is invertible, and therefore the matrix $\Gamma$ (defined by equation \ref{EQ-def-matrix-Gamma}) is invertible. Lemma \ref{lemma-Key-condition} implies that the entries of $\Gamma$ satisfy equation \ref{Key-EQ-coro}.

On the other hand if $\Gamma$ is an invertible matrix, and satisfies equation \ref{Key-EQ-coro}, then  by Lemma \ref{lemma-Key-condition}, it defines a morphism $\bgamma:\calA(T)\rightarrow \calA(S).$ Let $\Theta=[\theta_{ij}]$ be the inverse matrix of $\Gamma.$ Let us define the function $\theta:\calA(S)\rightarrow \calA(T)$ as follows
\begin{equation}\label{Eq-inverse-of-gamma}
\theta\left(\sum_{\boldsymbol{a}}\lambda_{\boldsymbol{a}} Y_1^{a_1}\cdots Y_n^{a_n}\right):= \sum_{\boldsymbol{a}}\lambda_{\boldsymbol{a}}\theta (Y_1)^{a_1}\cdots \theta (Y_n)^{a_n}
\end{equation}
where the sum runs over all the $\boldsymbol{a}=(a_1,\cdots,a_n)$ with $a_j\in\{0,1\},$ $\lambda_{\boldsymbol{a}}\in \F$ is any scalar and $\theta(Y_j)$ is given by the equation
\begin{equation}
\theta(Y_j):=\sum_{i=1}^{n}\theta_{ij}X_i.
\end{equation}
Is not difficult to see that $\theta$ is the inverse function of $\gamma.$ Since $\gamma$ is in particular an homomorphism of algebras, we conclude directly that $\theta$ is also an homomorphism of algebras. Since $\gamma$ preserves degrees, so does $\theta,$ therefore $\theta$ is a morphism in the category $\bcalNT$ and this implies that $\gamma$ is an isomorphism in the category $\bcalNT.$
\end{proof}

\subsection{Elementary triangular operations}\label{ssec-ETO}

In this section we define three matrix operations,  that we call \emph{elementary triangular operations} (ETO's). Basically each ETO transform a \nSLTM, $T$, into a new \nSLTM, $S.$ As we will see in Definitions \ref{def-ETO-scalar}, \ref{def-ETO-permuting} and \ref{def-ETO-Q}, each ETO has certain restrictions, that is, we only are allow to apply them on matrices $T,$ that satisfies a series of conditions.  Anyway, we shall prove that, whenever this transformation is well-defined, we obtain a matrix $S$ such that $T\sim S.$ Moreover we shall prove that if there is a finite sequence of, well-defined, ETO's that transform $T$ into $S,$ then we have that $T\sim S.$


\begin{definition}\label{def-ETO-scalar}
  Let $T=[t_{rk}]$ be a \nSLTM, fix an $1\leq r_1\leq n$ and an scalar $\alpha\neq 0.$ We define the matrix $S=\calP_{r_1}(T,\alpha),$ as the \nSLTM, whose entries $s_{rk}$ are given by:
  \begin{equation}\label{eq1-def-ETO-scalar}
    s_{rk}=\left\{\begin{array}{cc}
                    {t_{rk}}/{\alpha} & \textrm{if}\quad r=r_1 \\
                    \quad &\quad\\
                    \alpha t_{rk} & \textrm{if}\quad k=r_1 \\
                    \quad&\quad\\
                    t_{rk} & \textrm{otherwise}
                  \end{array}\right.
  \end{equation}
\end{definition}
\begin{example}
  Consider the matrix
  \begin{equation*}
    T=\left[\begin{matrix}
              0 & 0 & 0 & 0 & 0 & 0 \\
              1 & 0 & 0 & 0 & 0 & 0 \\
              2 &  1 & 0 & 0 & 0 & 0 \\
              3 & 3 & 1 & 0 & 0 & 0 \\
              1 & 3 & 2 & 3 & 0 & 0 \\
              2 & 2 & 3 & 2 &  1 & 0
            \end{matrix}\right]
  \end{equation*}
  then
  \begin{equation*}
    \calP_{3}(T,2)=\left[\begin{matrix}
              0 & 0 & 0 & 0 & 0 & 0 \\
              1 & 0 & 0 & 0 & 0 & 0 \\
              {\color{blue}1} & {\color{blue}{1}/{2}} & 0 & 0 & 0 & 0 \\
              3 & 3 & {\color{blue}2} & 0 & 0 & 0 \\
              1 & 3 & {\color{blue}4} & 3 & 0 & 0 \\
              2 & 2 & {\color{blue}6} & 2 &  1 & 0
            \end{matrix}\right],\quad \calP_{1}(T,3)=\left[\begin{matrix}
              0 & 0 & 0 & 0 & 0 & 0 \\
              {\color{red}3} & 0 & 0 & 0 & 0 & 0 \\
              {\color{red}6} &  1 & 0 & 0 & 0 & 0 \\
              {\color{red}9} & 3 & 1 & 0 & 0 & 0 \\
              {\color{red}3} & 3 & 2 & 3 & 0 & 0 \\
              {\color{red}6} & 2 & 3 & 2 &  1 & 0
            \end{matrix}\right]
  \end{equation*}
  \end{example}

\begin{lemma}\label{Lemma-ETO-P}
    Let $T=[t_{rk}]$ be a \nSLTM, fix $1\leq r_1\leq n$ and an scalar $\alpha\neq 0.$ Then $T\sim \calP_{r_1}(T,\alpha).$
\end{lemma}

\begin{proof}
  Let $S=\calP_{r_1}(T,\alpha)$ (as in definition \ref{def-ETO-scalar}).  Let $X_1,\dots,X_n$ and $Y_1,\dots,Y_n$ the generators of $\calA(T)$ and $\calA(S)$ respectively. We claim that the assignment
  \begin{equation}\label{eq1-lemma-ETO-escalar}
    X_r\mapsto \left\{\begin{array}{cc}
                        \alpha Y_r & \textrm{if}\quad r=r_1. \\
                        \quad&\quad\\
                        Y_r & \textrm{otherwise}
                      \end{array}\right.
  \end{equation}
  defines an isomorphism $\gamma:\calA(T)\rightarrow \calA(S).$ To see that, we only have to check that
\begin{equation}\label{eq2-lemma-ETO-escalar}
\gamma(X_r)^2=\sum_{j<r}{t_{rj}}\gamma(X_j)\gamma(X_r),
\end{equation}
    for each $1\leq r\leq n,$ and that the associated matrix $\Gamma$ is invertible. Now note that
  \begin{enumerate}
    \item If $r<r_1,$ then we have on one side
    \begin{equation*}
    \gamma(X_r)^{2}=Y_r^{2}=\sum_{j<r}s_{rj}Y_jY_r=\sum_{j<r}t_{rj}Y_jY_r,
    \end{equation*}
    the last equation is because $s_{rj}=t_{rj},$ whenever $r<r_1$ (see equation \ref{eq1-def-ETO-scalar}).

    On the other side, by equation  \ref{eq1-lemma-ETO-escalar}, we have:
    \begin{equation*}
    \sum_{j<r}{t_{rj}}\gamma(X_j)\gamma(X_r)=\sum_{j<r}{t_{rj}}Y_jY_r,
    \end{equation*}
    then it is clear that equation \ref{eq2-lemma-ETO-escalar} holds for this case.
    \item If $r=r_1,$ then we have on one side
    \begin{equation*}
      \gamma(X_{r_1})^2=(\alpha Y_r)^2=\alpha^2\sum_{j<r_1}{s_{r_1j}Y_jY_{r_1}}
      =\sum_{j<r_1}{\alpha t_{r_1j}Y_jY_{r_1}},
    \end{equation*}
    since $s_{r_{1}j}=t_{r_1j}/\alpha$ (see equation \ref{eq1-def-ETO-scalar}).

    On the other side, we have
    \begin{equation*}
      \sum_{j<r_1}{t_{r_1j}\gamma(X_{j})\gamma(X_{r_1})}
      =\sum_{j<r_1}{\alpha t_{r_1j}Y_jY_{r_1}}.
    \end{equation*}
    since $\gamma(X_j)=Y_j$ for $j<r_1$ and $\gamma(X_{r_1})=\alpha Y_{r_{1}}$ (see equation \ref{eq1-lemma-ETO-escalar}). Therefore equation \ref{eq2-lemma-ETO-escalar} holds for this case.
    \item if $r>r_1$ then we have on one side
    \begin{equation*}
      \gamma(X_r)^{2}=Y_r^2=\sum_{j\neq r_1}{s_{rj}Y_jY_r}+s_{rr_1}Y_{r_1}Y_{r}
      =\sum_{j\neq r_1}{t_{rj}Y_jY_r}+\alpha t_{rr_1}Y_{r_1}Y_{r},
    \end{equation*}
    since $s_{rj}=t_{rj},$ whenever $j\neq r_1$ and  $s_{rr_1}=\alpha t_{rr_1}$ (see equation \ref{eq1-def-ETO-scalar}).

    On the other side, we have
    \begin{equation*}
      \sum_{j<r}{t_{rj}\gamma(X_j)\gamma(X_r)}=\sum_{j\neq r_1}{t_{rj}Y_jY_r}+\alpha t_{rr_1}Y_{r_1}Y_{r},
    \end{equation*}
    since $\gamma(X_{j})=Y_j,$ whenever $j\neq r_1$ and $\gamma(X_{r_1})=\alpha Y_{r_1}.$ Therefore equation \ref{eq2-lemma-ETO-escalar} holds for this case.
  \end{enumerate}
  By the last analysis, we conclude that the assignment $\gamma$ defines an homomorphism $\calA(T)\rightarrow\calA(S).$ Is not difficult to see that the associated matrix $\Gamma=[\gamma_{ij}]$ is invertible, in fact, the entries of $\Gamma$ are given by

 \begin{equation}\label{eq00-ETO-matrices}
  \gamma_{ij}=\left\{\begin{array}{cc}
                  \alpha & \textrm{if}\quad i=j=r \\
                  1 & \textrm{if}\quad i=j\neq r \\
                  0 & \textrm{otherwise}
                \end{array}\right.
\end{equation}
  that is, $\Gamma=\alpha I_n,$  then $\det(\Gamma)=\alpha\neq 0,$ and we conclude that $\gamma$ defines an isomorphism $\calA(T)\rightarrow\calA(S).$
\end{proof}

In the following, we shall denote by $P_r(\alpha),$ the associated matrix of the isomorphism $\gamma$ defined in equation \ref{eq1-def-ETO-scalar}. For example if $n=6, r=3$ and $\alpha=2,$ we have:

\begin{equation}\label{eq001-ETO-matrices}
  P_r(\alpha):=\left[\begin{matrix}
  1&0& 0 & 0 & 0& 0\\
   0&1& 0 & 0 & 0& 0\\
   0&0& {\color{blue}2} & 0 & 0& 0\\
    0&0& 0 & 1 & 0& 0\\
     0&0& 0 & 0 & 1& 0\\
      0&0& 0 & 0 & 0& 1\\
  \end{matrix}\right]
  \end{equation}

Note that $P_r(\alpha)$ is an elementary matrix.

\begin{definition}\label{def-ETO-permuting}
  Let $T=[t_{rk}]$ be a \nSLTM, and fix $1\leq r_1<r_2\leq n.$ If the matrix $T$ satisfies the following conditions:
  \begin{enumerate}
    \item The entry $t_{r_2j}=0$ for each $r_1\leq j$
    \item The entry $t_{rr_1}=0$ for each $r_1< r< r_2$
  \end{enumerate}
  then we define the matrix $S=\calF_{(r_1,r_2)}(T),$ as the \nSLTM, whose entries $s_{rk}$ are given by:
  \begin{equation}\label{eq1-def-ETO-permuting}
    s_{rk}=\left\{\begin{array}{cc}
                    t_{r_2k} & \textrm{if}\quad r=r_1 \\
                    t_{r_1k} & \textrm{if}\quad r=r_2  \\
                    t_{rr_2} & \textrm{if}\quad k=r_1 \\
                    t_{rr_1} & \textrm{if}\quad k=r_2 \\
                    t_{rk} & \textrm{otherwise}
                  \end{array}\right.
  \end{equation}
\end{definition}
\begin{example}
  Consider the matrix
  \begin{equation*}
    T=\left[\begin{matrix}
              0 & 0 & 0 & 0 & 0 & 0 \\
              1 & 0 & 0 & 0 & 0 & 0 \\
              {\color{blue}2} &  {\color{blue}1} & 0 & 0 & 0 & 0 \\
              3 & 3 & 0 & 0 & 0 & 0 \\
              {\color{red}1} & {\color{red}3} & 0 & 0 & 0 & 0 \\
              2 & 2 & {\color{blue}3} & 2 &  {\color{red}1} & 0
            \end{matrix}\right]
  \end{equation*}
  and take $r_1=3$ and $r_2=5.$ Is not difficult to see that $T$ satisfies the conditions of Definition \ref{def-ETO-permuting}. Then the matrix $\calF_{(3,5)}(T)$ is well-defined. Moreover we have
  \begin{equation*}
    \calF_{(3,5)}(T)=\left[\begin{matrix}
              0 & 0 & 0 & 0 & 0 & 0 \\
              1 & 0 & 0 & 0 & 0 & 0 \\
              {\color{red}1} & {\color{red}3} & 0 & 0 & 0 & 0 \\
              3 & 3 & 0 & 0 & 0 & 0 \\
               {\color{blue}2} &  {\color{blue}1} & 0 & 0 & 0 & 0 \\
              2 & 2 &  {\color{red}1} & 2 & {\color{blue}3} & 0
            \end{matrix}\right]
  \end{equation*}
\end{example}

\begin{lemma}\label{Lemma-ETO-F}
  Let $T$ be a \nSLTM, let $1\leq r_1<r_2\leq n.$ If $T$ satisfies the conditions of Definition \ref{def-ETO-permuting}, then $T\sim\calF_{(r_1,r_2)}(T).$
\end{lemma}
\begin{proof}
  Let $S=\calF_{(r_1,r_2)}(T)$ (as in definition \ref{def-ETO-permuting}). Let $X_1,\dots,X_n$ and $Y_1,\dots,Y_n$ the generators of $\calA(T)$ and $\calA(S)$ respectively. We claim that the assignment
  \begin{equation}\label{eq1-lemma-ETO-F}
    X_r\mapsto\left\{\begin{array}{cc}
                       Y_{r_2} & \textrm{if}\quad r=r_1 \\
                       Y_{r_1} & \textrm{if}\quad r=r_2  \\
                       Y_r & \textrm{otherwise}
                     \end{array}\right.
  \end{equation}
  defines an isomorphism $\gamma:\calA(T)\rightarrow \calA(S).$ To see that, we only have to check that
\begin{equation}\label{eq2-lemma-ETO-F}
\gamma(X_r)^2=\sum_{j<r}{t_{rj}}\gamma(X_j)\gamma(X_r),
\end{equation}
    for each $1\leq r\leq n,$ and that the associated matrix $\Gamma$ is invertible. Now note that
  \begin{enumerate}
    \item For $r=r_1$ we have on one side
    \begin{equation*}
    \gamma(X_{r_1})^2=Y_{r_2}^2=\sum_{j<r_2}{s_{r_2j}Y_jY_{r_2}}
    =\sum_{j<r_1}{t_{r_1j}Y_jY_{r_2}}.
  \end{equation*}
  since $s_{r_2j}=t_{r_1j}$ (see equation \ref{eq1-def-ETO-permuting}), in particular $s_{r_2j}=t_{r_1j}=0$ for each $j\geq r_1,$ since $T$ is a \SLTM.  On the other side, by equation \ref{eq1-lemma-ETO-F}, we have
  \begin{equation*}
    \sum_{j<r_1}{t_{r_1j}\gamma(X_j)\gamma(X_{r_1})}
    =\sum_{j<r_1}{t_{r_1j}Y_jY_{r_2}}.
  \end{equation*}
  therefore equation \ref{eq2-lemma-ETO-F} holds in this case.
    \item For $r=r_2$ we have on one side
    \begin{equation*}
    \gamma(X_{r_2})^2=Y_{r_1}^2=\sum_{j<r_1}{s_{r_1j}Y_jY_{r_1}}
    =\sum_{j<r_1}{t_{r_2j}Y_jY_{r_1}}.
  \end{equation*}
  since $s_{r_1j}=t_{r_2j}$ (see equation \ref{eq1-def-ETO-permuting}), in particular,  $s_{r_1j}=t_{r_2j}=0$ for each $j\geq r_1,$ since the matrix $T$ has to satisfy the conditions on Definition \ref{def-ETO-permuting}.  On the other hand equation \ref{eq1-lemma-ETO-F} implies that
  \begin{equation*}
    \sum_{j<r_2}{t_{r_2j}\gamma(X_j)\gamma(X_{r_2})}
    =\sum_{j<r_1}{t_{r_2j}Y_jY_{r_1}}.
  \end{equation*}
  therefore, equation \ref{eq2-lemma-ETO-F} holds in this case.
    \item For $r<r_1$ it is immediate that equation \ref{eq2-lemma-ETO-F} holds, since by equation \ref{eq1-def-ETO-permuting}, $s_{rj}=t_{rj}$ for all $j$ and by equation \ref{eq1-lemma-ETO-F}, we have that $\gamma(X_j)=Y_j$ for all $j\geq r.$
    \item For $r_1<r<r_2$ we have on one side

     \begin{equation*}
    \gamma(X_r)^2=Y_r^2=\sum_{j<r}s_{rj}Y_jY_r=
    \sum_{j<r_1}{t_{rj}Y_jY_r}+\sum_{r_1<j<r}{t_{rj}Y_jY_r}
    \end{equation*}

the last equation is because, by equation \ref{eq1-def-ETO-permuting}, $s_{rj}=t_{rj}$ for all $j,$ in particular $s_{rr_1}=t_{rr_1}=0,$ since $T$ satisfies conditions on Definition \ref{def-ETO-permuting}. On the other side, the same argument and equation \ref{eq1-lemma-ETO-F} imply that
    \begin{equation*}
      \sum_{j<r}{t_{rj}\gamma(X_j)\gamma(X_r)}=\sum_{j<r_1}{t_{rj}Y_jY_r}+\sum_{r_1<j<r}{t_{rj}Y_jY_r}.
    \end{equation*}
    Therefore, equation \ref{eq2-lemma-ETO-F} holds in this case.
    \item For $r_2<r$ we have on one side:
    \begin{equation*}
    \begin{array}{cc}
      \gamma(X_r)^2 & =\sum_{j\neq r_1,r_2}{s_{rj}Y_jY_r}
    +s_{rr_1}Y_{r_1}Y_r+s_{rr_2}Y_{r_2}Y_r\\
      \quad & \quad \\
      \quad & =\sum_{j\neq r_1,r_2}{t_{rj}Y_jY_r}+t_{rr_2}Y_{r_1}Y_r+t_{rr_1}Y_{r_2}Y_r
    \end{array}
    \end{equation*}
    since  by equation \ref{eq1-def-ETO-permuting} $s_{rj}=t_{rj}$ for all $j\neq r_1,r_2,$ $s_{rr_1}=t_{rr_2}$ and $s_{rr_2}=t_{rr_1}.$ On the other side, we have:
    \begin{equation*}
      \sum_{j<r}{t_{rj}\gamma(X_j)\gamma(X_r)} =\sum_{j\neq r_1,r_2}{t_{rj}Y_jY_r} +t_{rr_1}Y_{r_2}Y_r
    +t_{rr_2}Y_{r_1}Y_r
    \end{equation*}
    since, by equation \ref{eq1-lemma-ETO-F},  $\gamma(X_j)=Y_j$ for all $j\neq r_1,r_2,$ $\gamma(X_{r_1})=Y_{r_2}$ and $\gamma(X_{r_2})=Y_{r_1}.$ Therefore equation \ref{eq2-lemma-ETO-F} holds in this case.
  \end{enumerate}
  The last analysis implies that the assignment $\gamma$ defines a homomorphism $\calA(T)\rightarrow \calA(S).$ Is not difficult to see that the associated matrix $\Gamma$ is invertible, in fact, the entries of $\Gamma$ are given by

  \begin{equation}\label{eq2-ETO-matrices}
     \gamma_{ij}=\left\{\begin{array}{cc}
                          1 & \textrm{if}\quad i=j\neq r_1,r_2 \\
                          1 & \textrm{if}\quad i=r_1,j=r_2 \\
                          1 & \textrm{if}\quad i=r_2,j=r_1 \\
                          0 & \textrm{otherwise}
                        \end{array}\right.
\end{equation}
  Therefore $\det(\Gamma)=-1$ and we conclude that $\gamma$ defines an isomorphism $\calA(T)\rightarrow \calA(S).$
\end{proof}

In the following we shall denote by $F_{(r_1,r_2)},$ the associated matrix of the isomorphism $\gamma$ defined in equation \ref{eq1-lemma-ETO-F}. For example if $n=6, r_1=2$ and $r_2=5,$ we have:

\begin{equation}\label{eq002-ETO-matrices}
  F_{(2,5)}:=\left[\begin{matrix}
  1&0& 0 & 0 & 0& 0\\
   0&0& 0 & 0 & {\color{blue}1}& 0\\
   0&0& {1} & 0 & {0}& 0\\
    0&0& 0 & 1 & 0& 0\\
     0& {\color{blue}1}&0 & 0 & 0& 0\\
      0&0& 0 & 0 & 0& 1\\
  \end{matrix}\right]
  \end{equation}

Note that $F_{(r_1,r_2)}$  is an elementary matrix.

Recall that for a given \nSLTM, $T=[t_{ij}]$ we denote $\Delta^{(1)}_{i,j,k}(T)=t_{ki}+t_{kj}t_{ji}$ (see subsection \ref{ssec-notations}).

\begin{definition}\label{def-ETO-Q}
Let $T=[t_{rk}]$ be a \nSLTM, take an scalar $\beta\in\F$ and assume that there are indices $1\leq k_0<r_0\leq n,$ such that $t_{r_0,j}=0$ for each $j>k_0.$ If one of the following conditions holds:
  \begin{enumerate}
    \item $k_0=1$
    \item $k_0>1$ and $\Delta^{(1)}_{i,k_0,r_0}(T)=\beta t_{k_0,i},$ for each $1\leq i<k_0$
  \end{enumerate}
  Then we define the matrix $S=\calQ_{(r_0,k_0)}(T,\beta)$ as the \nSLTM, whose engtries $s_{rk}$ are given by:
  \begin{equation}\label{eq1-def-ETO-Q}
    s_{rk}=\left\{\begin{array}{cc}
                    t_{r_0,k_0}-2\beta & \textrm{if}\quad (r,k)=(r_0,k_0) \\
                    \quad & \quad\\
                    t_{r,k_{0}}+\beta t_{r,r_{0}} & \textrm{if}\quad r>r_0\quad \textrm{and}\quad k=k_0  \\
                    \quad & \quad\\
                    t_{rk} & \textrm{otherwise}
                  \end{array}\right.
  \end{equation}
\end{definition}

\begin{example}
  Consider the matrix
  \begin{equation*}
    T=\left[\begin{matrix}
              0 & 0 & 0 & 0 & 0 \\
              2 & 0 & 0 & 0 & 0 \\
              1 & 0 & 0 & 0 & 0 \\
              3 & 1 & 2 & 0 & 0\\
              1 & 2 & 3 & 1 & 0
            \end{matrix}\right].
  \end{equation*}
  Let $(r_0,k_0)=(3,1)$ and $\beta=2.$ Note that the conditions of Definition \ref{def-ETO-Q} hold. In this case $k_0=1.$ Therefore we have
  \begin{equation*}
    \calQ_{(3,1)}(T,2)=\left[\begin{matrix}
              0 & 0 & 0 & 0 & 0 \\
              2 & 0 & 0 & 0 & 0 \\
              {\color{blue}1-4} & 0 & 0 & 0 & 0 \\
              {\color{blue}3+4} & 1 & 2 & 0 & 0\\
              {\color{blue}1+6} & 2 & 3 & 1 & 0
            \end{matrix}\right]=
            \left[\begin{matrix}
              0 & 0 & 0 & 0 & 0 \\
              2 & 0 & 0 & 0 & 0 \\
              {\color{blue}-3} & 0 & 0 & 0 & 0 \\
              {\color{blue}7} & 1 & 2 & 0 & 0\\
              {\color{blue}7} & 2 & 3 & 1 & 0
            \end{matrix}\right]
  \end{equation*}
\end{example}

\begin{example}
  Consider the matrix
  \begin{equation*}
    T=\left[\begin{matrix}
              0 & 0 & 0 & 0 & 0 \\
              2 & 0 & 0 & 0 & 0 \\
              1 & 2 & 0 & 0 & 0 \\
              3 & 6 & 2 & 0 & 0\\
              1 & 2 & 3 & 1 & 0
            \end{matrix}\right].
  \end{equation*}
  Let $(r_0,k_0)=(4,3)$ and $\beta=5.$ Note that the conditions of Definition \ref{def-ETO-Q} hold. In this case we have:
  \begin{equation*}
    \Delta^{(1)}_{1,3,4}=5=\beta t_{3,1},\quad \Delta^{(1)}_{2,3,4}=10=\beta t_{3,2}.
  \end{equation*}

   Therefore we have
  \begin{equation*}
    \calQ_{(4,3)}(T,5)=
            \left[\begin{matrix}
             0 & 0 & 0 & 0 & 0 \\
              2 & 0 & 0 & 0 & 0 \\
              1 & 2 & 0 & 0 & 0 \\
              3 & 6 & {\color{blue}2-10} & 0 & 0\\
              1 & 2 & {\color{blue}3+5} & 1 & 0
            \end{matrix}\right]=
            \left[\begin{matrix}
             0 & 0 & 0 & 0 & 0 \\
              2 & 0 & 0 & 0 & 0 \\
              1 & 2 & 0 & 0 & 0 \\
              3 & 6 & {\color{blue}-8} & 0 & 0\\
              1 & 2 & {\color{blue}8} & 1 & 0
            \end{matrix}\right]
  \end{equation*}
\end{example}

\begin{lemma}\label{lemma-ETO-Q}
Let $T=[t_{rk}]$ a \nSLTM, let $1\leq k_0<r_0\leq n$ be a pair of indices and $\beta\in \F$ an scalar. Assuming that the conditions of Definition \ref{def-ETO-Q} hold, then we have that $T\sim \calQ_{(r_0,k_0)}(T,\beta)$
\end{lemma}

\begin{proof}
We shall assume that $k_0>1.$ The case $k_0=1$ is easier and it is left to the readers.

 Assuming that $k_0>1$ and $\Delta^{(1)}_{i,k_0,r_0}=\beta t_{k_0,i},$ for each $1\leq i<k_0,$ then the matrix $S=\calQ_{(r_0,k_0)}(T,\beta)$ is well-defined and its entries $s_{rk}$ are given by equation \ref{eq1-def-ETO-Q} (see Definition \ref{def-ETO-Q}).  Let $X_1,\dots,X_n$ and $Y_1,\dots,Y_n$ the generators of $\calA(T)$ and $\calA(S)$ respectively.  We claim that the assignment:
  \begin{equation}\label{eq1-lemma-ETO-Q}
    X_r\mapsto \gamma(X_r)=\left\{\begin{array}{cc}
                        \beta Y_{k_0}+Y_{r_0} & \textrm{if}\quad r=r_0  \\
                        Y_{r} & \textrm{otherwise}
                      \end{array}\right.
  \end{equation}
  defines an isomorphism $\gamma:\calA(T)\rightarrow \calA(S).$ As in Lemmas \ref{Lemma-ETO-P} and \ref{Lemma-ETO-F}, we only  have to check that
\begin{equation}\label{eq2-lemma-ETO-Q}
\gamma(X_r)^2=\sum_{j<r}{t_{rj}}\gamma(X_j)\gamma(X_r),
\end{equation}
    for each $1\leq r\leq n,$ and that the associated matrix $\Gamma$ is invertible. Now note that:

   \begin{enumerate}
   \item If $1\leq r< r_0$ then it is clear that
  \begin{equation*}
    \gamma(X_r)^2=Y^2_{r}=\sum_{j<r}{s_{rj}Y_jY_r}=\sum_{j<r}{t_{rj}\gamma(X_j)\gamma(X_r)},
  \end{equation*}
  since by equation \ref{eq1-def-ETO-Q} $s_{rj}=t_{rj}$ and by equation \ref{eq1-lemma-ETO-Q}, $\gamma(X_j)=Y_j$ for all $j\leq r.$
  \item  If $r=r_0,$ we have on one side:

  \begin{equation*}
    \begin{array}{cc}
      \gamma(X_{r_0})^2 =\quad & (\beta Y_{k_0}+Y_{r_0})^2 \\
      \quad & \quad \\
      \quad\quad  =&\beta^2 Y_{k_0}^2+2\beta Y_{k_0}Y_{r_0}+Y_{r_0}^2\\
      \quad &\quad \\
      \quad\quad  =& \beta^2\sum_{j<k_0}{s_{k_0j}}Y_jY_{k_0}+\left(2\beta Y_{k_0}Y_{r_0}+\sum_{j<r_{0}}{s_{r_0j}Y_jY_{r_0}}\right) \\
      \quad & \quad \\
      \quad\quad  =& \beta^2\sum_{j<k_0}{t_{k_0j}}Y_jY_{k_0}+\sum_{j\leq k_0}{t_{r_0j}Y_jY_{r_0}}.
    \end{array}
  \end{equation*}
  The last step is because equation \ref{eq1-def-ETO-Q}, from where we know that $s_{k_0j}=t_{k_0j}$ for all $j,$ $s_{r_0j}=t_{r_0j}$ for $j\neq k_0,$ and  $s_{r_0k_0}=t_{r_0k_0}-2\beta.$  Also by the restrictions given in definition \ref{def-ETO-Q}, we have that $t_{r_0j}=0,$ for each $j>k_0.$

  On the other side, equation \ref{eq1-lemma-ETO-Q} implies that

  \begin{equation*}
    \begin{array}{ccc}
      \sum_{j<r_0}{t_{r_0j}\gamma(X_j)\gamma(X_{r_0})}& = & \sum_{j\leq k_0}{t_{r_0j}Y_j(\beta Y_{k_0}+ Y_{r_0})}\\
      \quad & \quad& \quad \\
      \quad & = & \beta(\sum_{j<k_0}{ t_{r_0j}Y_j Y_{k_0}+t_{r_0k_0} Y^2_{k_0})+\sum_{j\leq k_0}t_{r_0j}Y_jY_{r_0}}\\
      \quad & \quad& \quad \\
      \quad & = &\beta \sum_{j<k_0}( t_{r_0j}+t_{r_0k_0}t_{k_0j})Y_j Y_{k_0}+\sum_{j\leq k_0}t_{r_0j}Y_jY_{r_0}.
    \end{array}
  \end{equation*}

   Note that by our hypothesis, we have $\Delta^{(1)}_{j,k_0,r_0}(T)=t_{r_0j}+t_{r_0k_0}t_{k_0j}=\beta t_{k_0j},$ then we can reduce the expression on the parenthesis in the first sum to $\beta t_{k_0j},$ obtaining:
   \begin{equation*}
     \sum_{j<r_0}{t_{r_0j}\gamma(X_j)\gamma(X_{r_0})}=
     \beta^2\sum_{j<k_0}t_{k_0j}Y_j Y_{k_0}+\sum_{j\leq k_0}t_{r_0j}Y_jY_{r_0}
   \end{equation*}
  then equation \ref{eq2-lemma-ETO-Q} holds in this case.

 \item If $r>r_0$ we have on one side:

  \begin{equation*}
    \gamma(X_r)^2=Y_{r}^2=\sum_{j<r}{s_{rj}Y_jY_r}=(t_{rk_0}+\beta t_{rr_{0}})Y_{k_0}Y_r+\sum_{j\neq k_0}{t_{rj}Y_jY_r}
  \end{equation*}
  since by equation \ref{eq1-def-ETO-Q}, that $s_{rj}=t_{rj}$ for $j\neq k_0$ and $s_{rk_0}=t_{rk_0}+\beta t_{rr_0}.$

  On the other side, equation \ref{eq1-lemma-ETO-Q} implies that
  \begin{equation*}
  \begin{array}{cc}
  \sum_{j<r}{t_{rj}\gamma(X_j)\gamma(X_r)}=&t_{rk_0}Y_{k_0}Y_r+t_{rr_{0}}(\beta Y_{k_0}+Y_{r_{0}})Y_r+\sum_{j\neq k_0, r_0}{t_{rj}Y_jY_r}\\
  \quad &\quad\\
  \quad &= \beta^2\sum_{j<k_0}t_{k_0j}Y_j Y_{k_0}+\sum_{j\leq k_0}t_{r_0j}Y_jY_{r_0},
  \end{array}
  \end{equation*}
  therefore equation \ref{eq2-lemma-ETO-Q} holds in this case.
   \end{enumerate}

  By the above analysis we conclude that, the assignment $\gamma$ defines an homomorphism $\gamma:\calA(T)\rightarrow \calA(S).$ Is not difficult to see that the matrix $\Gamma$ is invertible, in fact, the entries of $\Gamma$ are given by:
\begin{equation}\label{eq3-ETO-matrices}
   \gamma_{ij}=\left\{\begin{array}{cc}
                          1 & \textrm{if}\quad i=j \\
                          \beta & \textrm{if}\quad i=k,j=r \\
                          0 & \textrm{otherwise}
                        \end{array}\right.
\end{equation}
  therefore $\det(\Gamma)=1\neq0,$ and we conclude that $\gamma$ defines an isomorphism.
\end{proof}

In the following we shall denote by $Q_{(r_0,k_0)}(\beta),$ the associated matrix of the isomorphism $\gamma$ defined in equation \ref{eq1-lemma-ETO-F}. For example if $n=6, r_0=4,k_0=2$ and $\beta=5,$ we have:

\begin{equation}\label{eq002-ETO-matrices}
  Q_{(4,2)}(5):=\left[\begin{matrix}
  1&0& 0 & 0 & 0& 0\\
   0&1& 0 & {\color{blue}5}  & 0& 0\\
   0&0& {1} & 0& 0& 0\\
    0&0& 0 & {1} & 0& 0\\
     0&0& 0 & 0 & {1}& 0\\
      0&0& 0 & 0 & 0& 1\\
  \end{matrix}\right]
  \end{equation}

Note that $Q_{(r_0,k_0)}(\beta),$   is an elementary matrix.

\begin{definition}\label{def-ETO-assignments}
We call \emph{Elementary triangular operations} (ETO for short) any of the following matrix assignments (whenever they are well-defined):

\begin{equation}\label{eq-ETO-assignments}
  \begin{array}{c}
    \calP_r(\quad,\alpha):T\mapsto \calP_r(T,\alpha)  \\
    \quad  \\
     \calF_{(r_1,r_2)}(\quad):T\mapsto \calF_{(r_1,r_2)}(T)\\
    \quad  \\
    \calQ_{(r_0,k_0)}(\quad,\beta):T\mapsto \calQ_{(r_0,k_0)}(T,\beta)
  \end{array}
\end{equation}

Given two \nSLTM, $T$ and $S$ we say that they are \emph{equivalent by ETO} if there is a finite sequence of ETO's that transform (progressively) $T$ into $S.$ We denote by $T\approx S$ whenever $T$ and $S$ are equivalent by ETO.
\end{definition}

It is clear that $\approx$ defines an equivalence relation on the set of all \nSLTM. In particular note that

\begin{equation}\label{eq-ETO-assignments-inv}
  \begin{array}{c}
    \calP_r(\quad,\alpha^{-1}):\calP_r(T,\alpha)\mapsto T  \\
    \quad  \\
     \calF_{(r_1,r_2)}(\quad):\calF_{(r_1,r_2)}(T)\mapsto T\\
    \quad  \\
    \calQ_{(r_0,k_0)}(\quad,-\beta):\calQ_{(r_0,k_0)}(T,\beta)\mapsto T
  \end{array}
\end{equation}

Each of the assignments in equation  \ref{eq-ETO-assignments-inv} is well-defined if and only if the corresponding assignment in equation \ref{eq-ETO-assignments} is well-defined.


\begin{example}
   Consider the matrices
  \begin{equation*}
    T=\left[\begin{matrix}
              0 & 0 & 0 & 0 & 0 \\
              2 & 0 & 0 & 0 & 0 \\
              1 & 0 & 0 & 0 & 0 \\
              3 & 1 & 2 & 0 & 0\\
              1 & 2 & 3 & 1 & 0
            \end{matrix}\right],\quad
    S=\left[\begin{matrix}
              0 & 0 & 0 & 0 & 0 \\
              -3 & 0 & 0 & 0 & 0 \\
              \frac{2}{5} & 0 & 0 & 0 & 0 \\
              7 & 2 & 5 & 0 & 0\\
              7 & 3 & 10 & 1 & 0
            \end{matrix}\right]
  \end{equation*}
  then $T\approx S.$ In fact one can check that
  \begin{equation*}
    T\approx T_2=\calQ_{(3,1)}(T,2)\approx T_3=\calF_{(2,3)}(T_2)\approx S=\calP_3(T_3,5).
  \end{equation*}
\end{example}


\begin{theorem}\label{theo-sim-vs-approx}
  Let $T,S$ be two \nSLTM. If $T\approx S,$ then $T\sim S.$
\end{theorem}
\begin{proof}
  If we assume that $T\approx S,$ that is, there is a finite sequence of ETO's, that transform $T$ into $S,$ say $T=T_0\rightarrow T_1:=\calE_1(T_0)\rightarrow T_2:=\calE_2(T_2)\rightarrow\cdots \rightarrow S=T_m:=\calE_m(T_{m-1}),$ where each $\calE_j$ is one of the three ETO's described in Definition \ref{def-ETO-assignments}, and the corresponding matrix $\calE_j(T_{j-1})$ is well-defined, that is, $T_{j-1}$ respects the corresponding restrictions of the operation $\calE_{j}.$ Then Lemmas \ref{Lemma-ETO-P}, \ref{Lemma-ETO-F} and \ref{lemma-ETO-Q} imply that at each step $T_{j-1}\sim T_{j}.$ The result follows by transitivity of relation $\sim.$

\end{proof}

\begin{proposition}\label{prop1-sim-vs-approx}
Let $T,S$ be two \nSLTM. Assuming that there is an isomorphism $\gamma:\calA(T)\rightarrow \calA(S),$ such that, the associated matrix $\Gamma$ is an elementary matrix, then $T\approx S.$
\end{proposition}

\begin{proof}
We fix a \nSLTM, T, and we left the matrix $S$ variable. Let $\Gamma$ be an elementary matrix. We separate our analysis in the following cases:
\begin{description}
\item[Case 1] If $\Gamma=P_{r_0}(\alpha)$ for some $1\leq r_0\leq n$ and some $\alpha\in \F^{\times}$ (for notation, see just below of the proof of Lemma \ref{Lemma-ETO-P}).  If we assume that $\Gamma$ defines an isomorphism $\gamma:\calA(T)\rightarrow \calA(S),$ then the entries $s_{ij}$ of $S$ must satisfy equation \ref{Key-EQ-coro}, where $t_{ij}$ are replaced by the entries of our fixed matrix $T$ and $\gamma_{ij}$ are replaced by the entries of $\Gamma,$ which are given explicitly in equation \ref{eq00-ETO-matrices}. Using equation \ref{Key-EQ-coro} in this way, we have:
    \begin{enumerate}
      \item If we replace a triple $(i,k,r)$ such that  $1\leq i<k\leq n,$ $1\leq r\leq n,$ and $k=r= r_0$  in equation \ref{Key-EQ-coro}, then we obtain
           \begin{equation}\label{eqK1C1-prop-sim-vs-approx}
             \alpha^{2}s_{r_0i}=t_{r_0i}\gamma_{ii}\gamma_{r_0r_0}=t_{r_0i}\alpha
           \end{equation}
           since in this situation, equation \ref{eq00-ETO-matrices} implies that $\gamma_{ir}=0,$ $\gamma_{ii}=1,$ $\gamma_{kr}=\alpha,$ and $\gamma_{kj}=0$ for all $j<r.$ Now from equation \ref{eqK1C1-prop-sim-vs-approx}, we obtain $s_{r_0i}=t_{r_{0}i}/\alpha,$ for all $i<r_0.$
       \item If we replace a triple $(i,k,r)$ such that  $1\leq i<k\leq n,$ $1\leq r\leq n,$ $k=r> r_0$ and $i= r_0$  in equation \ref{Key-EQ-coro}, then we obtain
           \begin{equation}\label{eqK2C1-prop-sim-vs-approx}
             s_{rr_0}=t_{rr_0}\gamma_{r_0r_0}\gamma_{rr}=t_{rr_0}\alpha.
           \end{equation}
           since by equation \ref{eq00-ETO-matrices}, we have in this situation that $\gamma_{ir}=0,$ $\gamma_{ii}=\alpha$ $\gamma_{kr}=\gamma_{rr}=1$ and $\gamma_{kj}=0$ for all $j<k.$
       \item For any other replacement of a triple $(i,k,r)$ such that  $1\leq i<k\leq n,$ $1\leq r\leq n,$ we obtain equations implying that $s_{ri}=t_{ri}$ or equations of the form $0=0.$
    \end{enumerate}
    The above analysis implies that the entries $s_{ki}$ of $S$ are given by equation \ref{eq1-def-ETO-scalar}, and then $S=\calP_{r_0}(T,\alpha).$ Therefore $S\approx T$ in this case.
    \item[Case 2] If $\Gamma=F_{(r_1,r_2)}$ for some $1\leq r_1<r_2\leq n$ (for notation, see just below of the proof of Lemma \ref{Lemma-ETO-F}). We apply a similar analysis as in \textbf{Case 1}, we assume that $\Gamma$ defines an isomorphism $\gamma:\calA(T)\rightarrow \calA(S),$ then the entries $s_{ij}$ of $S$ must satisfy equation \ref{Key-EQ-coro}, where $t_{ij}$ are replaced by the entries of our fixed matrix $T$ and $\gamma_{ij}$ are replaced by the entries of $\Gamma,$ which are given explicitly in equation \ref{eq2-ETO-matrices}. Using equation \ref{Key-EQ-coro} in this way, we have:
        \begin{enumerate}
          \item If we replace a triple $(i,k,r)$ such that  $1\leq i<k\leq n,$ $1\leq r\leq n,$ $k=r=r_2$ and $i=r_1$ in equation \ref{Key-EQ-coro}, we obtain
              \begin{equation}\label{eqR3C2-prop-sim-vs-approx}
                    0=t_{r_2r_1}\gamma_{r_2r_1}\gamma_{r_1r_2}=t_{r_2r_1}.
                  \end{equation}
            Note that this is one of the restrictions for the matrix $T,$ assumed in Definition \ref{def-ETO-permuting}.
          \item  If we replace a triple $(i,k,r)$ such that  $1\leq i<k\leq n,$ $1\leq r\leq n,$ $i=r_1,$  $r=r_2$ and $r_1<k<r_2$ in equation \ref{Key-EQ-coro}, we obtain
              \begin{equation}\label{eqR4C2-prop-sim-vs-approx}
                    0=t_{r_2k}\gamma_{kk}\gamma_{r_1r_2}=t_{r_2k},\quad (r_1<k<r_2).
                  \end{equation}
          and this is another of the restrictions for the matrix $T,$ assumed in Definition \ref{def-ETO-permuting}.
          \item If we replace a triple $(i,k,r)$ such that  $1\leq i<k\leq n,$ $1\leq r\leq n,$ $r_1<i=r<r_2$ and $k=r_2$ in equation \ref{Key-EQ-coro}, we obtain
              \begin{equation}\label{eqR1C2-prop-sim-vs-approx}
                0=t_{rr_1}\gamma_{r_2r_1}\gamma_{rr}=t_{rr_1}
              \end{equation}
              Note that this is the last of the restrictions for the matrix $T,$ assumed in Definition \ref{def-ETO-permuting}.
              \item If we replace a triple $(i,k,r)$ such that  $1\leq i<k\leq n,$ $1\leq r\leq n,$ $i=r=r_1$ and $k=r_2$ in equation \ref{Key-EQ-coro}, we obtain
              \begin{equation}\label{eqK0C2-prop-sim-vs-approx}
                s_{r_2r_1}=0
              \end{equation}
               In particular equations \ref{eqK0C2-prop-sim-vs-approx} imply that $s_{r_2r_1}=t_{r_1r_2}=0,$
          \item If we replace a triple $(i,k,r)$ such that  $1\leq i<k\leq n,$ $1\leq r\leq n,$ and $r_1<k=r<r_2$ and $i=r_1$ in equation \ref{Key-EQ-coro}, we obtain
                  \begin{equation}\label{eqR2C2-prop-sim-vs-approx}
                    s_{rr_1}=0,\quad (r_1<r<r_2).
                  \end{equation}
                In particular equations \ref{eqR1C2-prop-sim-vs-approx} and \ref{eqR2C2-prop-sim-vs-approx} imply that $s_{rr_1}=t_{rr_1},$ for all $r_1<r<r_2.$
          \item  If we replace a triple $(i,k,r)$ such that  $1\leq i<k\leq n,$ $1\leq r\leq n,$ and $r=r_1<i<k=r_2$ in equation \ref{Key-EQ-coro}, we obtain
              \begin{equation}\label{eqK5C2-prop-sim-vs-approx}
                    s_{r_2i}=0,\quad (r_1<i<r_2).
                  \end{equation}
                  In particular equation \ref{eqK5C2-prop-sim-vs-approx} implies that $s_{r_2i}=t_{r_1i}$ for all $(r_1<i<r_2).$
          \item If we replace a triple $(i,k,r)$ such that  $1\leq i<k\leq n,$ $1\leq r\leq n,$ $k=r_1$ and $r=r_2$ in equation \ref{Key-EQ-coro}, we obtain
              \begin{equation}\label{eqK4C2-prop-sim-vs-approx}
                    s_{r_1i}=t_{r_2i}\gamma_{ii}\gamma_{r_1r_2}=t_{r_2i}\quad (i<r_1).
                  \end{equation}
          \item  If we replace a triple $(i,k,r)$ such that  $1\leq i<k\leq n,$ $1\leq r\leq n,$ $i<r=r_1$ and $k=r_2$ in equation \ref{Key-EQ-coro}, we obtain
              \begin{equation}\label{eqK6C2-prop-sim-vs-approx}
                    s_{r_2i}=t_{r_1i}\gamma_{ii}\gamma_{r_2r_1}=t_{r_1i},\quad (i<r_1).
                  \end{equation}
          \item If we replace a triple $(i,k,r)$ such that  $1\leq i<k\leq n,$ $1\leq r\leq n,$ and $r_2<k=r$ and $i=r_1$ in equation \ref{Key-EQ-coro}, we obtain
              \begin{equation}\label{eqK2C2-prop-sim-vs-approx}
                    s_{rr_1}=t_{rr_2}\gamma_{r_1r_2}\gamma_{rr}=t_{rr_2}.
                  \end{equation}
          \item If we replace a triple $(i,k,r)$ such that  $1\leq i<k\leq n,$ $1\leq r\leq n,$ and $r_2<k=r$ and $i=r_2$ in equation \ref{Key-EQ-coro}, we obtain
              \begin{equation}\label{eqK2C2-prop-sim-vs-approx}
                    s_{rr_2}=t_{rr_1}\gamma_{r_1r_2}\gamma_{rr}=t_{rr_1}.
                  \end{equation}
          \item For any other replacement of a triple $(i,k,r)$ such that  $1\leq i<k\leq n,$ $1\leq r\leq n,$ we obtain equations implying that $s_{ri}=t_{ri}$ or equations of the form $0=0.$
        \end{enumerate}
        Note that equations \ref{eqR3C2-prop-sim-vs-approx}, \ref{eqR4C2-prop-sim-vs-approx} and \ref{eqR1C2-prop-sim-vs-approx} imply that if the matrix  $\Gamma=F_{(r_1,r_2)}$ defines an isomorphism $\gamma:\calA(T)\rightarrow \calA(S),$ then necessarily, the matrix $T$ must satisfy all the restrictions described in Definition \ref{def-ETO-permuting}. The above analysis imply that $S=\calF_{(r_1,r_2)}(T)$ and therefore $S\approx T.$
    \item[Case 3] If $\Gamma=Q_{(r_0,k_0)}(\beta),$  (for notation see just below of the proof of Lemma \ref{lemma-ETO-Q}). The analysis is similar to those of \textbf{Case 1} and \textbf{Case 2}, we have to replace the entries $\gamma_{ij}$ (defined in equation \ref{eq3-ETO-matrices}), in equation \ref{Key-EQ-coro}, collecting restrictions on $T$ and defining equations for the entries of $S.$ We shall assume that $k_0>1,$ the case $k_0=1$ is easier and it is left to the readers.

        In this case we have the following situations:
        \begin{enumerate}
        \item If we replace a triple $(i,k,r)$ such that $i=k_0,$  $r=r_0$ and $k_0<k<r_0$ in equation \ref{Key-EQ-coro}, we obtain
              \begin{equation}\label{eqR1C3-prop-sim-vs-approx}
                0=t_{r_0k}\gamma_{kk}\gamma_{k_0r_0}=t_{r_0k}\beta
              \end{equation}
              from where we obtain $t_{r_0k}=0$ for all $k>k_0.$ Note that this is one of the restrictions described in Definition \ref{def-ETO-Q}.
         \item If we replace a triple $(i,k,r)$ such that $i<k_0,$ and $k=r=k_0$  in equation \ref{Key-EQ-coro}, we obtain
              \begin{equation}\label{eqR0C3-prop-sim-vs-approx}
                s_{k_0i}=t_{k_0i}\gamma_{ii}\gamma_{k_0k_0}=t_{r_0i},\quad (1\leq i<k_0).
              \end{equation}
         \item If we replace a triple $(i,k,r)$ such that $i<k_0,$ $k=k_0$ and $r=r_0$ in equation \ref{Key-EQ-coro}, we obtain
             \begin{equation}\label{eqR2C3-prop-sim-vs-approx}
               \beta^{2}s_{k_0i}=t_{r_0k_0}s_{k_0i}\beta+t_{r_0i}\beta.
             \end{equation}
             Combining equations \ref{eqR0C3-prop-sim-vs-approx} and \ref{eqR2C3-prop-sim-vs-approx}, we obtain
             \begin{equation}\label{eqR3C3-prop-sim-vs-approx}
               \beta t_{k_0i}=\Delta_{i,k_0,r_0}^{(1)}(T),\quad (1\leq i<k_0)
             \end{equation}
             and this is another restriction described in Definition \ref{def-ETO-Q}.
          \item If we replace the triple $(i,k,r)=(k_0,r_0,r_0)$ in equation \ref{Key-EQ-coro}, we obtain
              \begin{equation}\label{eqK1C3-prop-sim-vs-approx}
                2\beta+s_{r_0k_0}=t_{r_0k_0}\gamma_{k_0k_0}\gamma_{r_0r_0}=t_{r_0k_0}
              \end{equation}
              from where we obtain $s_{r_{0}k_0}=t_{r_0k_0}-2\beta.$
          \item If we replace a triple $(i,k,r)$ such that $i=k_0,$  $k=r>r_0,$  in equation \ref{Key-EQ-coro}, we obtain
              \begin{equation}\label{eqK2C3-prop-sim-vs-approx}
                s_{rk_0}=t_{r_0k_0}\gamma_{k_0k_0}\gamma_{rr}
                +t_{rr_0}\gamma_{k_0r_0}\gamma_{rr}=t_{r_0k_0}+\beta t_{rr_0}.
              \end{equation}
         \item For any other replacement of a triple $(i,k,r)$ such that  $1\leq i<k\leq n,$ $1\leq r\leq n,$ we obtain equations implying that $s_{ri}=t_{ri}$ or equations of the form $0=0.$
        \end{enumerate}
        Note that equations \ref{eqR1C3-prop-sim-vs-approx}, and \ref{eqR3C3-prop-sim-vs-approx} imply that if the matrix  $\Gamma=Q_{(r_0,k_0)}(\beta)$ defines an isomorphism $\gamma:\calA(T)\rightarrow \calA(S),$ then necessarily, the matrix $T$ must satisfy all the restrictions described in Definition \ref{def-ETO-permuting}. The above analysis imply that $S=\calQ_{(r_0,k_0)}(T,\beta)$ and therefore $S\approx T.$
    \item[Case 4] If $\Gamma\neq P_{r_0}(\alpha),F_{(r_1,r_2)},Q_{(r_0,k_0)}(\beta),$ then necessarily $\Gamma=\left(Q_{(r_0,k_0)}(\beta)\right)^{t},$ that is, there is a pair $1\leq k_0< r_0\leq n$ and an scalar $\beta\neq 0$ such that
        \begin{equation}\label{eq-def-ETO-Q-trans}
          \gamma_{rk}=\left\{\begin{array}{cc}
                               1 & \textrm{if}\quad r=k \\
                               \beta & \textrm{if}\quad r=r_0,k=k_0 \\
                               0 & \textrm{otherwise}
                             \end{array}\right.
        \end{equation}
        Working exactly as we did in the three previous cases, we obtain that:
         \begin{enumerate}
           \item The entries of $T$ have the following restrictions:
           \begin{equation}\label{K1C4-prop-sim-vs-approx}
             \left\{\begin{array}{cc}
                      t_{rk_0}=0 & \textrm{for all}\quad k_0<r<r_0 \\
                      \quad & \quad\\
                      t_{r_0k}=0 & \textrm{for all}\quad k_0<k<r_0 \\
                      \quad & \quad\\
                      t_{k_0j}=-\beta t_{r_0j} & \textrm{for all}\quad 1\leq j<k_0\\
                      \quad & \quad\\
                       t_{r_0k_0}=2/\beta & \quad
                    \end{array}\right.
           \end{equation}
           In particular, they imply that $\Delta_{j,k_0,r_0}^{(1)}(T)=\frac{1}{\beta}t_{k_0j},$ for all $1\leq j<k_0,$ and therefore the matrix $T'=\calQ_{(r_0,k_0)}(T,1/\beta)$ is well-defined.
           \item The entries of $S$ depends on $T$ as follows:
           \begin{equation}\label{K1C4-prop-sim-vs-approx}
             s_{rk}=\left\{\begin{array}{cc}
                             -2/\beta & \textrm{if}\quad r=r_0,k=k_0\\
                              \quad & \quad\\
                             t_{k_0k}=-\beta t_{r_0k} & \textrm{if}\quad r=k_0,k<r_0 \\
                             \quad & \quad\\
                             -t_{r_0k} & \textrm{if}\quad r=r_0,k<r_0 \\
                             \quad & \quad\\
                              \beta t_{rk_0}+t_{rr_0}& \textrm{if}\quad r>r_0,k=r_0\\
                              \quad & \quad\\
                               t_{rk_0}& \textrm{if}\quad r>r_0,k=k_0\\
                              \quad & \quad\\
                             t_{rk} & \textrm{otherwise}
                           \end{array}\right.
           \end{equation}
         \end{enumerate}
         Now note that the entries $t_{rk}'$ of the matrix $T'=\calQ_{(r_0,k_0)}(T,\beta)$ are given by:
         \begin{equation}\label{K1C4-prop-sim-vs-approx}
           t_{rk}'=\left\{\begin{array}{cc}
                             0 & \textrm{if}\quad r=r_0,k=k_0\\
                              \quad & \quad\\
                             t_{r_0k} & \textrm{if}\quad r=r_0,k<r_0 \\
                             \quad & \quad\\
                             -\beta t_{r_0k} & \textrm{if}\quad r=k_0,k<r_0 \\
                             \quad & \quad\\
                              t_{rk_0}+t_{rr_0}/\beta& \textrm{if}\quad r>r_0,k=k_0\\
                              \quad & \quad\\
                             t_{rk} & \textrm{otherwise}
                           \end{array}\right.
         \end{equation}
         In particular we have that $t_{r_0k}'=t_{r_0k}=0$ for all $k>k_0$ and  $t'_{rk_0}=t_{rk_0}=0$ for all $k_0<r\leq r_0.$ This implies that the matrix $U:=\calF_{(k_0,r_0)}(T')$ is well-defined.

         Now note that the entries of $U$ are given by:
         \begin{equation}\label{K2C4-prop-sim-vs-approx}
           u_{rk}=\left\{\begin{array}{cc}
                             0 & \textrm{if}\quad r=r_0,k=k_0\\
                              \quad & \quad\\
                             t_{r_0k} & \textrm{if}\quad r=k_0,k<r_0 \\
                             \quad & \quad\\
                             -\beta t_{r_0k} & \textrm{if}\quad r=r_0,k<r_0 \\
                             \quad & \quad\\
                              t_{rk_0}+ t_{rr_0}/\beta& \textrm{if}\quad r>r_0,k=r_0\\
                              \quad & \quad\\
                              t_{rr_0}& \textrm{if}\quad r>r_0,k=k_0\\
                              \quad & \quad\\
                             t_{rk} & \textrm{otherwise}
                           \end{array}\right.
         \end{equation}
         Let us define the matrix $U'=\calP_{r_0}(\beta).$ Since the $\calP-$operation have no restrictions, $U'$ is well-defined and its entries are given by:
         \begin{equation}\label{K3C4-prop-sim-vs-approx}
           u_{rk}'=\left\{\begin{array}{cc}
                             0 & \textrm{if}\quad r=r_0,k=k_0\\
                              \quad & \quad\\
                             t_{r_0k} & \textrm{if}\quad r=k_0,k<r_0 \\
                             \quad & \quad\\
                             -t_{r_0k} & \textrm{if}\quad r=r_0,k<r_0 \\
                             \quad & \quad\\
                              \beta t_{rk_0}+t_{rr_0}& \textrm{if}\quad r>r_0,k=r_0\\
                              \quad & \quad\\
                              t_{rr_0}& \textrm{if}\quad r>r_0,k=k_0\\
                              \quad & \quad\\
                             t_{rk} & \textrm{otherwise}
                           \end{array}\right.
         \end{equation}
         Note in particular that $\Delta_{k,k_0,r_0}^{(1)}(U')=-u_{k_0k},$ for all $1\leq k<k_0,$ therefore we can define the matrix $V:=\calQ_{(r_0,k_0)}(U,-1)$ whose entries are given by:
         \begin{equation}\label{K4C4-prop-sim-vs-approx}
           v_{rk}=\left\{\begin{array}{cc}
                             2 & \textrm{if}\quad r=r_0,k=k_0\\
                              \quad & \quad\\
                             t_{r_0k} & \textrm{if}\quad r=k_0,k<r_0 \\
                             \quad & \quad\\
                             -t_{r_0k} & \textrm{if}\quad r=r_0,k<r_0 \\
                             \quad & \quad\\
                              \beta t_{rk_0}+t_{rr_0}& \textrm{if}\quad r>r_0,k=r_0\\
                              \quad & \quad\\
                              -\beta t_{rk_0}& \textrm{if}\quad r>r_0,k=k_0\\
                              \quad & \quad\\
                             t_{rk} & \textrm{otherwise}
                           \end{array}\right.
         \end{equation}
         Finally we define the matrix $S'=\calP_{k_0}(-1/\beta).$ By contruction $T\approx S'$ and one can check that the entries of $S'$ are given by
         \begin{equation}\label{K5C4-prop-sim-vs-approx}
           s_{rk}'=\left\{\begin{array}{cc}
                             -2/\beta & \textrm{if}\quad r=r_0,k=k_0\\
                              \quad & \quad\\
                             t_{k_0k}=-\beta t_{r_0k} & \textrm{if}\quad r=k_0,k<r_0 \\
                             \quad & \quad\\
                             -t_{r_0k} & \textrm{if}\quad r=r_0,k<r_0 \\
                             \quad & \quad\\
                              \beta t_{rk_0}+t_{rr_0}& \textrm{if}\quad r>r_0,k=r_0\\
                              \quad & \quad\\
                               t_{rk_0}& \textrm{if}\quad r>r_0,k=k_0\\
                              \quad & \quad\\
                             t_{rk} & \textrm{otherwise}
                           \end{array}\right.
         \end{equation}
         Comparing equations \ref{K1C4-prop-sim-vs-approx} and \ref{K5C4-prop-sim-vs-approx}, one conclude that $S=S'$ and therefore $T\approx S$ as desired.
\end{description}
\end{proof}

Note that in \textbf{Case 4} of the proof of Proposition \ref{prop1-sim-vs-approx}, the matrix $\Gamma=(\calQ_{(r_0,k_0)}(\beta))^{t}$ can be factored as follows:
\begin{equation}\label{eq1-Fact-Gamma-elementary}
  \Gamma=P_{k_0}(1/\beta)\cdot Q_{(r_0,k_0)}(-1)\cdot P_{r_0}(\beta)\cdot F_{(r_0,k_0)}\cdot Q_{(r_0,k_0)}(1/\beta)
\end{equation}
and that factorization yields the following \emph{well-defined} process of transformations:

\begin{equation}\label{eq2-Fact-Gamma-elementary}
\begin{array}{c}
  T\rightarrow   T':=\calQ_{(r_0,k_0)}(T,1/\beta)\rightarrow U:=\calF_{(k_0,r_0)}(T')
  \rightarrow \cdots \\
  \quad     \\
  \quad  \cdots \rightarrow U':=\calP_{r_0}(U,\beta)\rightarrow   V:=\calQ_{(r_0,k_0)}(U',-1)\rightarrow \cdots  \\
  \quad \quad \\
  \quad  \quad \quad  \cdots\rightarrow S=\calP_{k_0}(V,1/\beta).\\
\end{array}
\end{equation}

That is, the restrictions and relations for matrices $T,S$ obtained, when the entries of $\Gamma$ are replaced in the system of equations \ref{Key-EQ-coro}, are coherent with the process of transformations induced by the factorization given in equation \ref{eq1-Fact-Gamma-elementary}.

More generally, any invertible matrix $\Gamma$ can be factored as a product of elementary matrices, and each elementary matrix can be written as product of matrices of the forms $P_{r_0}(\alpha),$ $F_{(r_1,r_2)}$ or $ Q_{(r_0,k_0)}.$ Now if we assume that $\Gamma$ is an invertible matrix, that define an isomorphism $\gamma:\calA(T)\rightarrow \calA(S),$ for certain variable matrices $T,S,$ then the system of equations \ref{Key-EQ-coro}, will imply certain restrictions and relations on $T$ and $S.$ Based in the proof of Proposition \ref{prop1-sim-vs-approx}, it is sensible to think that, there will be a factorization $\Gamma=E_1\cdots E_m$ with matrices $E_j$ of one of the forms  $P_{r_0}(\alpha),$ $F_{(r_1,r_2)}$ or $ Q_{(r_0,k_0)},$ yielding a \emph{well-defined} process of transformations:
\begin{equation*}
  T\rightarrow T_2=\calE_1(T)\rightarrow\cdots \rightarrow S=\calE_m(T_{m-1}).
\end{equation*}

At this point we do not have a well proof of that fact, but in an enormous set of explicit examples, we always obtain evidences of its veracity. Based on that, we define the following:

\begin{conjecture}\label{conj-sim-vs-approx}\textbf{(Isomorphism vs Equivalence by ETO's)}
Let $T,S$ be two \nSLTM. If $T\sim S,$ then $T\approx S.$
\end{conjecture}

In the following section we shall provide more evidences of the veracity of Conjecture \ref{conj-sim-vs-approx}, see Corollaries \ref{coro-Prop-Nn-for-n2}, \ref{coro-Prop-Nn-for-n3} and Proposition \ref{prop2-sim-vs-approx}.

\section{Applications}\label{sec-int-applic}

In this section we apply the isomorphism criteria obtained in Section \ref{sec-Isom-criteria}, to give some first steps in the search of a complete classification of the objects of the category $\bcalNT.$

\subsection{Small cases}\label{ssec-small-cases}

In the smallest case, that is, $n=2$ we can easily count exactly the number of classes:
\begin{proposition}\label{Prop-Nn-for-n2}
  There is only one class in $\mathcal{TM}_2$ under relation $\sim.$
\end{proposition}
\begin{proof}
Let
\begin{equation}\label{eq1-Prop-Nn-for-n2}
  T=\left[\begin{matrix}
                  0 & 0 \\
                  t & 0
                \end{matrix}\right],\quad S=\left[\begin{matrix}
                                                    0 & 0 \\
                                                    s & 0
                                                  \end{matrix}\right].
\end{equation}
  Then the matrix
  \begin{equation}\label{eq2-Prop-Nn-for-n2}
  \Gamma=\left[\begin{matrix}
                                   1 & \frac{t-s}{2} \\
                                   0 & 1
                                 \end{matrix}\right],
\end{equation}
  defines an isomorphism $\calA(T)\rightarrow \calA(S).$
\end{proof}

\begin{corollary}\label{coro-Prop-Nn-for-n2}
  If $T,S$ are any two matrices in  $\mathcal{TM}_2.$  If $T\sim S$ then $T\approx S.$
\end{corollary}
\begin{proof}
  Let $T,S$ as in equation \ref{eq1-Prop-Nn-for-n2}. Proposition \ref{Prop-Nn-for-n2} says that always $T\sim S.$ Now equation  \ref{eq2-Prop-Nn-for-n2}, implies that $S=\calQ_{(2,1)}(T,\frac{t-s}{2}),$ and therefore $T\approx S.$
\end{proof}

The case $n=3$ is a little bit more complicated, but again, we can count exactly  the number of classes. For the rest of this subsection for a given $3-$\SLTM, $U=[u_{ij}],$ we denote $\Delta=2u_{31}+u_{32}u_{21}=\Delta^{(2)}_{1,2,3}(U).$ (for notations, see subsection \ref{ssec-notations}).

\begin{proposition}\label{prop-Cl-eq-2-for-n3}
  There are exactly two classes in $\mathcal{TM}_3$ under relation $\sim.$
\end{proposition}

\begin{proof}
  We divide the proof in a series of steps:
  \begin{description}
    \item[Step 1] Note that the matrices
  $$0_3=\left[\begin{matrix}                 0 & 0 & 0 \\
                                           0 & 0 & 0 \\
                                           0 & 0 & 0
                                         \end{matrix}\right],\quad B_{3,2}=\left[\begin{matrix}
                                           0 & 0 & 0 \\
                                           0 & 0 & 0 \\
                                           1 & 1 & 0
                                         \end{matrix}\right]$$
  are not related by $\sim.$ In fact if there was an invertible $3\times3-$ matrix $\Gamma=[\gamma_{ij}]$ satisfying equation \ref{Key-EQ-coro}, then we have for any $1\leq r\leq 3$
  \begin{equation}\label{EQr-lemma-n3}
    \left\{\begin{array}{c}
             2\gamma_{1r}\gamma_{2r}=0 \\
             \gamma_{3r}(2\gamma_{1r}+\gamma_{3r})=0  \\
             \gamma_{3r}(2\gamma_{2r}+\gamma_{3r})=0
           \end{array}\right.
  \end{equation}
  If we assume that $\gamma_{3r}\neq 0,$ then by the second and third equations in \ref{EQr-lemma-n3}, we have   $\gamma_{1r}\neq 0$ and $\gamma_{2r}\neq 0,$ but this contradicts the  first equation in \ref{EQr-lemma-n3}, therefore $\gamma_{3r}=0$ for all $1\leq r \leq 3,$ implying that $\Gamma$ is not invertible.

  Note that \textbf{Step 1}, implies that there are at least two classes in in $\mathcal{TM}_3$ under relation $\sim.$
  \item[Step 2] Given a $3-$\SLTM, $U=[u_{ij}]$ such that $u_{32}\neq 0$ and $\Delta\neq 0,$ then we have   $U\sim B_{3,2}.$ In fact, if we consider the matrix:
\begin{equation*}
  \Gamma= \left[\begin{matrix}
                   2/\Delta & u_{21}/\Delta & 0 \\
                   0 & 1/u_{32} & 0 \\
                   0 & 0 & 1
                 \end{matrix}\right],
  \end{equation*}
one can check that the entries of $\Gamma$ satisfy equation \ref{Key-EQ-coro} of Theorem \ref{coro-Key-condition}, (taking $T=U$ and $S=B_{3,2}$). One can also check that $\det(\Gamma)\neq 0,$ therefore $\Gamma$ defines an isomorphism $\gamma:\calA(U)\rightarrow \calA(B_{3,2}).$
\item[Step 3] Given a $3-$\SLTM, $U=[u_{ij}]$ such that $u_{32}=0$ or $\Delta=0,$ then $U\sim 0_3.$ In fact, the following list shows all the different cases for $U\neq 0_3$ satisfying the hypothesis of this lemma, and the corresponding matrix $\Gamma$ that defines an isomorphism $\gamma:\calA(U)\rightarrow \calA(0_3):$
  \begin{equation*}
    \begin{array}{ccc}
      U=\left[\begin{matrix}
        0 & 0 & 0 \\
        u_{21} & 0 & 0 \\
        0 & 0 & 0
      \end{matrix}\right] & (u_{21}\neq 0,u_{32}=0), & \Gamma=\left[\begin{matrix}
        1 & {u_{21}}/2 & 0 \\
        0 & 1 & 0 \\
        0 & 0 & 1
      \end{matrix}\right] \\
      \quad & \quad & \quad  \\
      U=\left[\begin{matrix}
        0 & 0 & 0 \\
        0 & 0 & 0 \\
        u_{31} & 0 & 0
      \end{matrix}\right] & (u_{31}\neq 0,u_{32}=0), & \Gamma=\left[\begin{matrix}
        1 & 0 & u_{31}/2 \\
        0 & 1 & 0 \\
        0 & 0 & 1
      \end{matrix}\right] \\
      \quad & \quad & \quad  \\
      U=\left[\begin{matrix}
        0 & 0 & 0 \\
        0 & 0 & 0 \\
        0 & u_{32} & 0
      \end{matrix}\right] & (u_{32}\neq 0,\Delta=0), & \Gamma=\left[\begin{matrix}
        1 & 0 & 0 \\
        0 & 1 & u_{32}/2 \\
        0 & 0 & 1
      \end{matrix}\right] \\
      \quad & \quad & \quad  \\
      U=\left[\begin{matrix}
        0 & 0 & 0 \\
        u_{21} & 0 & 0 \\
        u_{31} & 0 & 0
      \end{matrix}\right] & (u_{21},u_{31}\neq 0,u_{32}=0), & \Gamma=\left[\begin{matrix}
        1 & u_{21}/2 & u_{31}/2 \\
        0 & 1 & 0 \\
        0 & 0 & 1
      \end{matrix}\right] \\
      \quad & \quad & \quad  \\
      U=\left[\begin{matrix}
        0 & 0 & 0 \\
        u_{21} & 0 & 0 \\
        u_{31} & u_{32} & 0
      \end{matrix}\right] & (u_{21},u_{31},u_{32}\neq 0,\Delta=0), & \Gamma=\left[\begin{matrix}
        1 & u_{21}/2 & 0 \\
        0 & 1 & u_{32}/2 \\
        0 & 0 & 1
      \end{matrix}\right] \\
    \end{array}
  \end{equation*}
  One can check in each case that $\Gamma$ is an invertible matrix and that its entries satisfy equation \ref{Key-EQ-coro} of Theorem \ref{coro-Key-condition}, (taking $T=U$ and $S=0_3$). Note that in the last case, we have $u_{31}=-\frac{1}{2}u_{32}u_{21},$ since $\Delta=0.$
  \item[Step 4] Let $U$ be a $3-$\SLTM. By  \textbf{Step 2} and \textbf{Step 3}, we have that $U$ is in the class of $0_3$ or in the class of $B_{3,2}.$ Finally \textbf{Step 1} implies that there are exactly two classes in $\mathcal{TM}_3.$
  \end{description}
\end{proof}


\begin{corollary}\label{coro-Prop-Nn-for-n3}
  If $T,S$ are any two matrices in  $\mathcal{TM}_3.$  If $T\sim S,$ then $T\approx S.$
\end{corollary}

\begin{proof}
  Let $U$ be a $3-$\SLTM. By the proof of Proposition \ref{prop-Cl-eq-2-for-n3}, we have that:
  \begin{itemize}
    \item If $U\sim B_{3,2},$ then one can check that the matrix $\Gamma$ defined in \textbf{Step 2} can be factored as follows:
        \begin{equation*}
          \Gamma=\left[\begin{matrix}
                         2/\Delta & 0 & 0 \\
                         0 & 1 & 0 \\
                         0 & 0 & 1
                       \end{matrix}\right]\cdot \left[\begin{matrix}
                         1 & 0 & 0 \\
                         0 & 1/u_{32} & 0 \\
                         0 & 0 & 1
                       \end{matrix}\right]
                       \cdot \left[\begin{matrix}
                         1 & u_{32}/2 & 0 \\
                         0 & 1 & 0 \\
                         0 & 0 & 1
                       \end{matrix}\right],
        \end{equation*}
        that is
        \begin{equation}\label{eq1-coro-Prop-Nn-for-n3}
          \Gamma=P_{1}(2/\Delta)\cdot P_{2}(1/u_{32})\cdot Q_{(2,1)}(u_{21}/2).
        \end{equation}
        Now this factorization yields the following well-defined process of transformations:
        \begin{equation}\label{eq2-coro-Prop-Nn-for-n3}
          U\rightarrow U_1:=\calQ_{(2,1)}(u_{21}/2)\rightarrow
          U_2:=\calP_{2}(U_1,1/u_{32})\rightarrow U_3=\calP_{1}(U_2,2/\Delta).
        \end{equation}
        Is not difficult to check that $U_3=B_{3,2},$ and therefore one conclude that $U\approx B_{3,2}.$
    \item If $U\sim 0_{3},$ then one can check that in all the cases described in \textbf{Step 3}, one have $U\approx 0_3.$ In fact all the corresponding matrices $\Gamma,$ can be factored in a product of matrices of the form $Q_{(r_0,k_0)}(\beta)$ and the factorization yields a well-defined process of transformations, via $\calQ-$operations, $U\rightarrow \cdots \rightarrow 0_n,$ therefore $U\sim 0_n$ implies that $U\approx 0.$ We left details to the readers.
  \end{itemize}
  Let $T,S\in \mathcal{TM}_3$ such that $T\sim S.$ As a consequence of the proof of Proposition \ref{prop-Cl-eq-2-for-n3}, this implies that one and only one of the following is true $T\sim S\sim 0_3,$ or $T\sim S\sim B_{3,2}.$ Finally the above analysis implies in both cases that $T\approx S$ whenever $T\sim S.$
\end{proof}

\subsection{The class of $0_n.$ }\label{ssec-class-of-0n}
We would like to understand, how the information coming from the codifying matrices, can be used to describe completely an specific isomorphism class in $\bcalNT.$ In this direction, in this section we describe a very explicit criterion (Theorem \ref{theo-zero-class-n-gen}) for determining when a given \nSLTM, say $U,$ is in the same class of the matrix $0=0_n.$ This criterion is a generalization of \textbf{Step 3} in the proof of Proposition \ref{prop-Cl-eq-2-for-n3}.

Along this subsection we denote $X_1,\cdots X_n$ and $Y_1,\cdots Y_n$ the generators of $\calA(0_n)$ and $\calA(U),$ respectively. Also recall from subsection \ref{ssec-notations}, the notation $\Delta_{i,j,k}^{(2)}(U)$ (or simply $\Delta^{(2)}_{i,j,k}$) given by $\Delta_{i,j,k}^{(2)}(U)=2u_{ki}+u_{kj}u_{ji}.$





\medskip


\begin{lemma}\label{lemma-zero-class-necessary-condition}
  If $U=[u_{ij}]$ is a \nSLTM, such that $0_n\sim U,$ then one of the following conditions holds
   \begin{description}
    \item[1] $u_{kj}=0,$ for each pair of indices $1< j<k\leq n.$
    \end{description}
    or
    \begin{description}
    \item[2] For each pair of indices $1< j<k\leq n$ such that $u_{kj}\neq 0,$ we have that $\Delta^{(2)}_{i,j,k}=0,$ for each  $1\leq i<j.$
    \end{description}
\end{lemma}

\begin{proof}
  If $\gamma:\calA(0_n)\rightarrow\calA(U),$ is an isomorphism, then by Theorem \ref{coro-Key-condition} we have
  \begin{equation}\label{Key-EQ-adapted1}
    2\gamma_{ir}\gamma_{kr}+\gamma^{2}_{kr}u_{ki}
    =\gamma_{kr}(2\gamma_{ir}+\gamma_{kr}u_{ki})=0
  \end{equation}
  for each $1\leq r\leq n,$ and each $1\leq i<k\leq n.$ Let us assume that there is a $u_{kj}\neq 0,$ for some pair $(k,j)$ with $1<j<k.$ Since the matrix $\Gamma=[\gamma_{ij}]$ is invertible, there is a $r=1,\dots,n$ such that  $\gamma_{kr}\neq 0,$ therefore equation \ref{Key-EQ-adapted1} implies that
  \begin{equation*}
    u_{ki}=-2\frac{\gamma_{ir}}{\gamma_{kr}},\quad \textrm{for each}\quad 1\leq i\leq j.
  \end{equation*}
  On the other hand if we take $i=j$ in equation \ref{Key-EQ-adapted1}, we can also conclude that $\gamma_{jr}\neq 0,$ therefore
  \begin{equation}\label{Key-EQ-adapted2}
    \gamma_{jr}(2\gamma_{ir}+\gamma_{jr}u_{ji})=0,\quad (1\leq i<j)
  \end{equation}
  since equation \ref{Key-EQ-adapted1} is true for each for each $1\leq r\leq n,$ and each $1\leq i<k\leq n.$ Now equation \ref{Key-EQ-adapted2} implies that
  \begin{equation*}
    u_{ji}=-2\frac{\gamma_{ir}}{\gamma_{jr}},\quad \textrm{for each}\quad 1\leq i< j.
  \end{equation*}
  Therefore we have
  \begin{equation*}
    \Delta^{(2)}_{i,j,k}=2u_{ki}+u_{kj}u_{ji}
    =-4\frac{\gamma_{ir}}{\gamma_{kr}}
    +4\frac{\gamma_{jr}}{\gamma_{kr}}\frac{\gamma_{ir}}{\gamma_{jr}}=0,\quad \textrm{for each}\quad 1\leq i< j.
  \end{equation*}
\end{proof}

\begin{definition}\label{def-leaders-chains-etc}
  Let $U\neq 0_n,$ be a \nSLTM.
  \begin{enumerate}
    \item A \emph{leader} of $U$ is a pair $(k,j)$ such that $u_{kj}\neq 0$ and $u_{kc}=0$ for all $c>j.$ We denote by $\calL_U$ the set of all leaders of $U.$
    \item We define the \emph{graph} of leaders of $U$ as the directed graph $\Omega_{U}=(V_U,A_U),$ where the set of vertices is $V_U=\{1,\dots,n\}$ and where arrows are defined by:
        \begin{equation*}
          k\rightarrow j\quad \textrm{if and only if}\quad (k,j)\in\calL_{U}.
        \end{equation*}
        We also denote by $k\twoheadrightarrow j$ if there is a finite sequence $k_1,\dots,k_m$ such that $k_1=k,$ $k_m=j$ and $k_i\rightarrow k_{i+1}$ for each $i=1,\dots,m-1.$
    \item Let $\Theta=(V,A)$ a connected subgraph of $\Omega_{U}.$ We denote by $\calC_{\Theta}$ the set given by:
        \begin{equation*}
          \calC_{\Theta}=\{(k,j):k,j\in V\quad\textrm{and}\quad k\rightarrow j\}.
        \end{equation*}
        If the set $\calC_{\Theta}$ is not empty, we say that it is a \emph{chain of leaders} of $U.$ If $\Theta$ is a (non-trivial) connected component of $\Omega_{U},$ then we say that $\calC_{\Theta}$ is a \emph{maximal} chain of leaders. If $\Theta_1,\Theta_2$ are two different connected components of $\Omega,$ then we say that $\calC_{\Theta_1}$ and $\calC_{\Theta_2}$ are two independent chains of leaders of $U.$

    \item A \emph{simple chain of leaders} of $U$ is a subset $\calC\subset\calL_{U}$ say $\calC=\{(k_1,j_1),(k_2,j_2),\dots,(k_l,j_l)\}$ with $l\geq 1$ and $1\leq k_1<k_2<\cdots<k_l\leq n$ such that
        \begin{equation*}
          k_{r}= j_{r+1}\quad \textrm{for each}\quad  1\leq r\leq l.
        \end{equation*}
        \item A \emph{ramified chain of leaders} of $U$ is any chain of leaders that is not simple.
  \end{enumerate}
\end{definition}

\begin{lemma}\label{lemma-properties-of-Omega}
Given a \nSLTM, $U,$ we have that:
\begin{enumerate}
  \item\label{it1-lemma-properties-of-Omega} Relation $k\rightarrow j$ satisfy the following properties:
   \begin{equation}\label{eq-propiedad-sumidero0}
   k\rightarrow j\Rightarrow k>j.
   \end{equation}
  and
   \begin{equation}\label{eq-propiedad-sumidero}
   (k\rightarrow j_1\wedge k\rightarrow j_2)\Rightarrow j_1=j_2.
   \end{equation}

  \item\label{it2-lemma-properties-of-Omega} The following relation:
    \begin{equation}\label{eq-def-unrhd}
    k\unrhd j\quad\textrm{if and only if}\quad (k=j\vee k\twoheadrightarrow j)
    \end{equation}
    is a partial order on the set of vertices $V_U$ of $\Omega_{U}.$
  \item\label{it3-lemma-properties-of-Omega} For each $k\in V_U$ there is a unique minimal element $k_0$ of $V_U$ (relative to $\unrhd$), such that $k\unrhd k_0.$ Moreover if $k_1,k_2\in V_U$ are in the same connected component of $\Omega_{U},$ then there is a unique minimal element $k_0$ of $V_U$ such that $k_1\unrhd k_0$ and $k_2\unrhd k_0.$
  \item\label{it4-lemma-properties-of-Omega} If $k$ in not a minimal element of $V_U$ under $\unrhd,$ then the set $$\calC_{k}=\{(r,c):(k=r\vee k\twoheadrightarrow r)\wedge r\rightarrow c \}$$ is a simple chain of leaders of $U.$
  \item $k\in V_{U}$ is a minimal element if and only if there is no $j<k$ such that $(k,j)\in\calL_{U}$ and that means that the $k-$th row of $U$ is zero.
\end{enumerate}
\end{lemma}

The proof of Lemma \ref{lemma-properties-of-Omega} is direct and it is left to the readers.


\begin{example}\label{ex-mega-ex-chains}
   Consider the matrix:
  \begin{equation}\label{eq6-ex-mega-ex-chains}
    U=\left[\begin{array}{cccccccccccc}
              0 & 0 & 0 & 0 & 0 & 0 & 0 & 0 & 0 & 0 & 0 & 0 \\
              {\color{blue}\boldsymbol{1}} & 0 & 0 & 0 & 0 & 0 & 0 & 0 & 0 & 0 & 0 & 0 \\
              -1 & {\color{blue}\boldsymbol{2}} & 0 & 0 & 0 & 0 & 0 & 0 & 0 & 0 & 0 & 0 \\
              0 & 0 & 0 & 0 & 0 & 0 & 0 & 0 & 0 & 0 & 0 & 0  \\
              1 & -2 & {\color{blue}2} & 0 & 0 & 0 & 0 & 0 & 0 & 0 & 0 & 0  \\
              \frac{1}{2} & -1 & {\color{violet}1} & 0 & 0 & 0 & 0 & 0 & 0 & 0 & 0 & 0  \\
              -1 & 2 & -2 & 0 & {\color{blue}2} & 0 & 0 & 0 & 0 & 0 & 0 & 0 \\
              0 & 0 & 0 & {\color{red}1} & 0 & 0 & 0 & 0 & 0 & 0 & 0 & 0 \\
              0 & 0 & 0 & -\frac{3}{2} & 0 & 0 & 0 & {\color{red}3} & 0 & 0 & 0 & 0 \\
              \frac{1}{2} & -1 & 1 & 0 & -1 & 0 & {\color{blue}1} & 0 & 0 & 0 & 0 & 0 \\
              0 & 0 & 0 & 0 & 0 & 0 & 0 & 0 & 0 & 0 & 0 & 0\\
              -\frac{3}{4} & \frac{3}{2} & -\frac{3}{2} & 0 & 0 & {\color{violet}3} & 0 & 0 & 0 & 0 & 0 & 0
            \end{array}\right]
  \end{equation}
  Then the set of leaders is $\calL_{U}=\{(2,1),(3,2),(5,3),(6,3),(7,5),(8,4),(9,8),(10,7),(12,6)\}.$ The directed graph associated is given by:
  
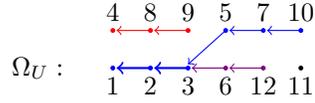
\begin{equation*}
    \begin{tikzpicture}[xscale=0.5,yscale=0.5]
\node[] at (-1,0) {$\Omega_{U}:$};
\node[below] at (1,0) {$1$};
\node[below] at (2,0) {$2$};
\node[below] at (3,0) {$3$};
\node[below] at (4,0) {$6$};
\node[below] at (5,0) {$12$};
\node[above] at (4,1) {$5$};
\node[above] at (5,1) {$7$};
\node[above] at (6,1) {$10$};
\node[above] at (1,1) {$4$};
\node[above] at (2,1) {$8$};
\node[above] at (3,1) {$9$};
\node[below] at (6,0) {$11$};
\draw[fill,blue] (1,0) circle [radius=0.05];
\draw[fill,red] (1,1) circle [radius=0.05];
\draw[fill,red] (2,1) circle [radius=0.05];
\draw[fill,red] (3,1) circle [radius=0.05];
\draw[fill] (6,0) circle [radius=0.05];
\draw[fill,blue] (2,0) circle [radius=0.05];
\draw[fill,blue] (3,0) circle [radius=0.05];
\draw[fill,violet] (4,0) circle [radius=0.05];
\draw[fill,violet] (5,0) circle [radius=0.05];
\draw[fill,blue] (4,1) circle [radius=0.05];
\draw[fill,blue] (5,1) circle [radius=0.05];
\draw[fill,blue] (6,1) circle [radius=0.05];
\draw[->,blue,thick] (2,0)--(1.1,0);
\draw[->,blue,thick] (3,0)--(2.1,0);
\draw[->,violet] (4,0)--(3.1,0);
\draw[->,violet] (5,0)--(4.1,0);
\draw[->,blue] (4,1)--(3,0.1);
\draw[->,blue] (5,1)--(4.1,1);
\draw[->,blue] (6,1)--(5.1,1);
\draw[->,red] (2,1)--(1.1,1);
\draw[->,red] (3,1)--(2.1,1);
    \end{tikzpicture}
  \end{equation*}
  Note that $\calL_{U}$ is divided in two independent maximal chains of leaders, $\calL_{U}=\calC_1\cup\calC_2,$ where
  $$\calC_1=\{(2,1),(3,2),(5,3),(6,3),(7,3),(10,7),(12,6)\},\quad \calC_2=\{(8,4),(9,8)\}.$$

  Also note that $\calC_2$ is a simple chain, while $\calC_1$ is a ramified chain, such that $\calC_1=\calC^{(1)}\cup \calC^{(2)},$ where $\calC^{(1)}, \calC^{(2)}$ are the simple chains of leaders given by: $$\calC^{(1)}=\{(2,1),(3,2),(5,3),(7,3),(10,7)\},\quad \calC^{(2)}=\{(2,1),(3,2),(6,3),(12,6)\}.$$
  Finally note that the minimal vertices of $\Omega_U$ are $1,4,11$, this is because the $1$st, the $4$th and the $11$th rows of $U$ are zero row. In particular the vertex $11$ is an isolated vertex. 
\end{example}

\begin{lemma}\label{lemma-chains-vs-Delta-zero}
  Let $U$ be a \nSLTM. If for each pair $(r,j)\in \calL_{U}$ with $j>1,$ we have that $\Delta_{i,j,r}^{(2)}=0,$ for each $i<j,$ then for each pair $(r,c)$ with $1<c<j$ such that $u_{rc}\neq 0,$ we have that $\Delta_{i,c,r}^{(2)}=0$ for each $i<c.$
\end{lemma}

\begin{proof}
  Let  $(r,j)\in \calL_{U}$ and $1<c<j$ such that $u_{rc}\neq 0.$  By our hypothesis, $\Delta_{i,j,r}^{(2)}=2u_{ri}+u_{rj}u_{ji}=0$ for each $i<j,$ and this implies that $u_{ri}=0$ if and only if $u_{ji}=0.$ In particular, we have that $u_{jc}\neq 0,$ since $u_{rc}\neq0.$

 Fix an $i<c.$  We have two cases:

  \begin{enumerate}
   \item If we assume that $(j,c)\in\calL_{U}$  (that is $r\rightarrow j\rightarrow c$), then by our hypothesis,  we have that:
       \begin{equation*}
         \Delta_{c,j,r}^{(2)}=2u_{rc}+u_{rj}u_{jc}=0,\quad
         \Delta_{i,c,j}^{(2)}=2u_{ji}+u_{jc}u_{ci}=0,\quad
          \Delta_{i,j,r}^{(2)}=2u_{ri}+u_{rj}u_{ji}=0
       \end{equation*}
      Then we have:
       \begin{equation*}
         u_{rc}=-\frac{1}{2}u_{rj}u_{jc},\quad
         u_{ji}=-\frac{1}{2}u_{jc}u_{ci},\quad
         2u_{ri}=\frac{1}{2}u_{rj}u_{jc}u_{ci},
       \end{equation*}
      therefore
      \begin{equation*}
        \Delta_{i,c,r}^{(2)}=2u_{ri}+u_{rc}u_{ci}
        =\frac{1}{2}u_{rj}u_{jc}u_{ci}-\frac{1}{2}u_{rj}u_{jc}u_{ci}
        =0.
      \end{equation*}
   \item If $(j,c)\notin\calL_{U},$ then since $u_{jc}\neq 0$ then necessarily there is a $c<j_2<j$ such that $(j,j_2)\in \calL_{U}$ and again this implies that $u_{j_2,c}\neq 0.$ By a recursive argument, we can find a finite sequence: $c=j_m<\cdots <j_2<j_1=j$ such that $(j_s,j_{s+1})\in \calL_{U}$ for each $s=1,\dots,m-1,$ that is $j=j_1\rightarrow j_2\rightarrow\cdots \rightarrow j_m=c,$ (in other words  $r\rightarrow j\twoheadrightarrow c$).

   In this case, by our hypothesis,  we have:
       \begin{equation}\label{eq2-1-proof-lemma-chains-vs-Delta-zero}
         \Delta_{i,j_1,r}^{(2)}=2u_{r,i}+u_{r,j_1}u_{j_1,i}=0,\quad
         \Delta_{j_m,j_1,r}^{(2)}=2u_{r,j_m}+u_{r,j_1}u_{j_1,j_m}=0,
       \end{equation}
       \begin{equation}\label{eq2-2-proof-lemma-chains-vs-Delta-zero}
         \Delta_{i,j_{s},j_{s-1}}^{(2)}=2u_{j_{s-1},i}+u_{j_{s-1},j_s}u_{j_s,i}=0,\quad
         \textrm{for each}\quad s=2,\dots,m.
       \end{equation}
       \begin{equation}\label{eq2-3-proof-lemma-chains-vs-Delta-zero}
         \Delta_{j_m,j_{s},j_{s-1}}^{(2)}=
         2u_{j_{s-1},j_m}+u_{j_{s-1},j_s}u_{j_s,j_m}=0,\quad
         \textrm{for each}\quad s=2,\dots,m-1.
       \end{equation}
       From the equations in \ref{eq2-1-proof-lemma-chains-vs-Delta-zero} we obtain respectively:
       \begin{equation}\label{eq2-4-proof-lemma-chains-vs-Delta-zero}
         2u_{r,i}=-u_{r,j_1}u_{j_1,i}\quad\textrm{and}\quad
         u_{r,j_m}=-\frac{1}{2}u_{r,j_1}u_{j_1,j_m}
       \end{equation}
       Using recursively equations \ref{eq2-2-proof-lemma-chains-vs-Delta-zero} and \ref{eq2-3-proof-lemma-chains-vs-Delta-zero} we obtain respectively:

       \begin{equation}\label{eq2-5-proof-lemma-chains-vs-Delta-zero}
         u_{j_1,i}=\left(-\frac{1}{2}\right)^{m-1}u_{j_1,j_2}u_{j_2,j_3}\cdots u_{j_{m},i}.
       \end{equation}
       and
       \begin{equation}\label{eq2-6-proof-lemma-chains-vs-Delta-zero}
         u_{j_1,j_m}=\left(-\frac{1}{2}\right)^{m-2}u_{j_1,j_2}u_{j_2,j_3}\cdots u_{j_{m-1},j_m}.
       \end{equation}
       therefore, if we replace respectively in equations \ref{eq2-5-proof-lemma-chains-vs-Delta-zero} and \ref{eq2-6-proof-lemma-chains-vs-Delta-zero}, the equations of \ref{eq2-4-proof-lemma-chains-vs-Delta-zero}, we obtain

        \begin{equation}\label{eq2-7-proof-lemma-chains-vs-Delta-zero}
         2u_{r,i}=-\left(-\frac{1}{2}\right)^{m-1}u_{r,j_1}u_{j_1,j_2}\cdots u_{j_{m-1},j_m}u_{j_m,i}.
       \end{equation}
       and
       \begin{equation}\label{eq2-8-proof-lemma-chains-vs-Delta-zero}
         u_{r,j_m}=\left(-\frac{1}{2}\right)^{m-1}u_{r,j_1}u_{j_1,j_2}\cdots u_{j_{m-1},j_m}.
       \end{equation}
       Finally, by equations \ref{eq2-7-proof-lemma-chains-vs-Delta-zero} and \ref{eq2-8-proof-lemma-chains-vs-Delta-zero}, we have:
       \begin{equation*}
         \Delta_{i,c,r}^{(2)}=\Delta_{i,j_m,r}^{(2)}=2u_{r,i}+u_{r,j_m}u_{j_m,i}=0
       \end{equation*}
       as desired.
  \end{enumerate}
\end{proof}

\begin{lemma}\label{lemma-zero-class-one-chain}
  Assuming that $\calL_{U}$ is itself a simple chain of leaders. If for each $(k,j)\in \calL_{U}$ and each $1<c<j$  we have $\Delta^{(2)}_{c,j,k}(U)=0,$ then $0_n\sim U.$
\end{lemma}
\begin{proof}
Recall that $X_1,\dots, X_n$ and $Y_1,\dots, Y_n$ denote the generators of $\calA(0_n)$ and $\calA(U),$ respectively. Let $\calL_{U}=\{(k_1,j_1),(k_2,j_2),\dots,(k_l,j_l)\},$ then the hypothesis implies that the presentation of $\calA(U)$ is given by the relations
  \begin{equation}\label{EQ-present-onechain-AU}
    Y_r^2=\left\{\begin{array}{cc}
                   \sum_{j\leq j_i}u_{rj}Y_jY_r & \textrm{if} \quad r=k_i\in\{k_1,\dots,k_l\} \\
                   \quad & \quad \\
                   0 & \textrm{otherwise}
                 \end{array}\right.
  \end{equation}
  Since $\calL_{U}$ is itself a simple chain, we have that $k_{i}=j_{i+1}$ for each $i=1,\dots,l-1.$ We also denote $k_0:=j_1.$ We define for each minimal element $r\in V_U$ (in particular for $r=k_0$),  the element $\gamma(X_r)=Y_r.$ Being minimal means that there is no $c<r$ such that $(r,c)\in \calL_{U},$ in particular by Lemma \ref{lemma-properties-of-Omega}, we know that $\gamma(X_r)^2=Y_r^{2}=0.$

  For $r\in\{k_1,\dots,k_l\}$ we define the element $\gamma(X_r)\in\calA(U),$ recursively by:

  \begin{equation}\label{EQ2-0-present-onechain-AU}
    \gamma(X_{k_i})=\left\{\begin{array}{cc}
                             -\frac{1}{2}u_{{k_1}{k_0}}Y_{k_0}+Y_{k_1} & \textrm{if}\quad i=1.  \\
                            \quad & \quad  \\
                             -\frac{1}{2}u_{{k_i}{k_{i-1}}}\gamma(X_{k_{i-1}})+Y_{k_i} & \textrm{if}\quad i>1.
                           \end{array}\right.
  \end{equation}
  We claim that the assignment $X_r\mapsto \gamma(X_r)$ defines an isomorphism $\gamma:\calA(0_n)\rightarrow \calA(U).$ By equation \ref{EQ-present-onechain-AU}, it is enough to prove that $\gamma(X_{k_i})^2=0$ for each $i=1,\dots,l.$

  For $i=1,$ note that $Y_{k_0}^2=0,$ and $Y_{k_1}^2=u_{{k_1}{k_0}}Y_{k_0}Y_{k_1},$ then
  \begin{equation*}
    \gamma(X_{k_1})^2=\left(-\frac{1}{2}u_{{k_1}{k_0}}Y_{k_0}+Y_{k_1}\right)^2
    =-u_{{k_1}{k_0}}Y_{k_0}Y_{k_1}+Y_{k_1}^2
    =-u_{{k_1}{k_0}}Y_{k_0}Y_{k_1}+u_{{k_1}{k_0}}Y_{k_0}Y_{k_1}=0.
  \end{equation*}
  For $i>1,$ we assume that $\gamma(X_{k_{i-1}})^2=0,$ then
  \begin{equation}\label{EQ2-present-onechain-AU}
    \gamma(X_{k_i})^2=
    -u_{{k_i}{k_{i-1}}}\gamma(X_{k_{i-1}})Y_{k_i}+Y_{k_i}^2
    =-u_{{k_i}{k_{i-1}}}\gamma(X_{k_{i-1}})Y_{k_i}+\sum_{j\leq j_i}u_{{{k_i}}j}Y_jY_{k_i}
  \end{equation}

Recall that we want to prove that $ \gamma(X_{k_i})^2=0,$ then we need to examine the two terms on the right of equation \ref{EQ2-present-onechain-AU} to see them clearer:

On one hand, one can check recursively (similarly to the proof of Lemma \ref{lemma-chains-vs-Delta-zero}), that:

  \begin{equation*}
   \gamma(X_{k_{i}})=\sum_{j=0}^{i}\left[\left(-\frac{1}{2}\right)^{i-j}\left(\prod_{s=j+1}^{i}u_{k_{s},k_{s-1}} \right)Y_{k_{j}}\right]
  \end{equation*}
  Therefore

  \begin{equation}\label{EQ2-2-present-onechain-AU}
  {\begin{array}{ccc}
  -u_{{k_i}{k_{i-1}}}\gamma(X_{k_{i-1}})Y_{k_i}&=&-u_{{k_i}{k_{i-1}}}
   \sum_{j=0}^{i-1}\left[\left(-\frac{1}{2}\right)^{i-1-j}\left(\prod_{s=j+1}^{i-1}u_{k_{s},k_{s-1}} \right)Y_{k_{j}}\right]Y_{k_i}\\
   \quad&\quad&\quad\\
   \quad&=&-\sum_{j=0}^{i-1}\left[\left(-\frac{1}{2}\right)^{i-1-j}\left(\prod_{s=j+1}^{i}u_{k_{s},k_{s-1}} \right)Y_{k_{j}}Y_{k_i}\right]
  \end{array}}
  \end{equation}

   On the other hand, note that condition $\Delta^{(2)}_{j,k_{i-1},k_i}=2u_{{k_i},j}+u_{{k_i},k_{i-1}}u_{k_{i-1},j}=0,$ for each $1\leq j<j_i=k_{i-1},$ implies that

   \begin{equation}\label{EQ2-3-present-onechain-AU}
     \sum_{j\leq j_i}u_{{{k_i}}j}Y_jY_{k_i}=\sum_{j=0}^{i-1}u_{k_{i},k_{j}}Y_{k_{j}}Y_{k_i}
   \end{equation}

   Moreover, we can use equation $2u_{{k_i},j}+u_{{k_i},k_{i-1}}u_{k_{i-1},j}=0,$ recursively to check that

   \begin{equation}\label{EQ3-present-onechain-AU}
     u_{k_{i},k_{j}}=\left(-\frac{1}{2}\right)^{i-1-j}\left(\prod_{s=j+1}^{i}u_{k_{s},k_{s-1}} \right),
   \end{equation}
   for any $k_0\leq k_{j}\leq k_{i-2}.$
   In fact, assuming that equation \ref{EQ3-present-onechain-AU} is true, for any $k_0\leq k_{j}\leq k_{i-2},$ and noticing that that $u_{k_{i},k_{i-2}}=\left(-\frac{1}{2}\right)u_{k_{i},k_{i-1}}u_{k_{i-1},k_{i-2}},$ we obtain:

   \begin{equation*}
   \begin{array}{ccc}
   u_{k_{i+1},k_{j}}&=&\left(-\frac{1}{2}\right)u_{k_{i+1},k_{i}}\left(-\frac{1}{2}\right)^{i-1-j}
     \left(\prod_{s=j+1}^{i}u_{k_{s},k_{s-1}} \right)\\
     \quad&\quad&\quad\\
     \quad&=&\left(-\frac{1}{2}\right)^{i-j}\left(\prod_{s=j+1}^{i+1}u_{k_{s},k_{s-1}} \right).
   \end{array}
   \end{equation*}

   Therefore equation \ref{EQ3-present-onechain-AU} is true for any for any $k_0\leq k_{j}\leq k_{i-1}.$

   Finally using equations \ref{EQ2-2-present-onechain-AU}, \ref{EQ2-3-present-onechain-AU}, \ref{EQ3-present-onechain-AU}, one can check that equation \ref{EQ2-present-onechain-AU} reduces to $\gamma(X_{k_i})^2=0,$ as desired. Now note that from equation \ref{EQ2-0-present-onechain-AU}, we have that $\det(\Gamma)=1$ and therefore $\Gamma$ defines an isomorphism $\gamma:\calA(0_n)\rightarrow \calA(U).$
\end{proof}

\begin{proposition}\label{prop-zero-class-general}
  Let $U\neq 0_n,$ be a \nSLTM, and assume that for each $(k,j)\in \calL_{U},$ with $j>1,$ we have that $\Delta_{i,j,k}^{(2)}=0,$ for each $i<j,$ then $0_n\sim U.$
\end{proposition}

\begin{proof}
  For each minimal element $k\in V_{U},$ we define $\gamma(X_k)=Y_k.$ Again, being minimal means that there is no a $j<k$ such that $(k,j)\in \calL_{U},$ therefore Lemma \ref{lemma-properties-of-Omega} implies that $\gamma(X_k)^{2}=Y_k^2=0.$

  For $k\in V_{U}$ not minimal, we define recursively:
  \begin{equation}\label{eq-recursive-chains-general}
    \gamma(X_k)=-\frac{u_{kj}}{2}\gamma(X_j)+Y_k,\quad\textrm{if}\quad k\rightarrow j.
  \end{equation}
  By equation \ref{eq-propiedad-sumidero} of Lemma \ref{lemma-properties-of-Omega}, $\gamma(X_k)$ is well-defined.

  By item \ref{it4-lemma-properties-of-Omega} of Lemma \ref{lemma-properties-of-Omega}, we know that  the set $\calC_{k}=\{(r,c):(k=r\vee k\twoheadrightarrow r)\wedge r\rightarrow c \}$ is a simple chain of leaders, therefore the fact that $\gamma(X_k)^{2}=0$ follows analogously to the proof of this in Lemma \ref{lemma-zero-class-one-chain}. Note that each independent chain can be treated independently. Therefore we have that the assignment $X_r\mapsto \gamma(X_r)$ defines an homomorphism $\gamma:\calA(0_n)\rightarrow \calA(U).$ It is clear that the associated matrix $\Gamma$ is invertible, therefore $\gamma:\calA(0_n)\rightarrow \calA(U)$ is an isomorphism.
\end{proof}

\begin{theorem}\label{theo-zero-class-n-gen}
  Given a \nSLTM, $U\neq 0_n,$ then $0_n\sim U$ if and only if one of the following conditions holds
  \begin{description}
    \item[1]\label{theo-zero-class-n-gen-It1} $u_{kj}=0,$ for each pair of indices $1< j<k\leq n.$
    \end{description}
    or
    \begin{description}
    \item[2] \label{theo-zero-class-n-gen-It2}  For each $(k,j)\in \calL_{U},$ with $j>1,$ we have that $\Delta_{i,j,k}^{(2)}=0,$ for each $i<j.$
    \end{description}
\end{theorem}


\begin{proof}
  Lemma \ref{lemma-zero-class-necessary-condition} show that if $0_n\sim U,$ then necessarily one of the conditions \textbf{1} or \textbf{2} holds.  Proposition \ref{prop-zero-class-general} together with Lemma \ref{lemma-chains-vs-Delta-zero} implies that if one of the conditions \textbf{1} or \textbf{2} holds, then necessarily $0_n\sim U.$
\end{proof}

\begin{example}\label{ex2-mega-ex-chains}
  Let $U$ be the matrix of equation \ref{eq6-ex-mega-ex-chains}. Note that $U$  satisfies the hypothesis of Proposition \ref{prop-zero-class-general}. As we have mentioned in Example \ref{ex-mega-ex-chains}, $\calL_{U}=\calC_1\cup\calC_2$ is divided in two independent chains, a ramified one $\calC_1=\calC^{(1)}\cup \calC^{(2)},$ with $\calC^{(1)}=\{(2,1),(3,2),(5,3),(7,3),(10,7)\},$ $\calC^{(2)}=\{(2,1),(3,2),(6,3),(12,6)\}$ and a simple one $\calC_2=\{(8,4),(9,8)\}.$

  In this case the minimal vertices of $\Omega_{U}$ are $1,4,11,$ where only $11$ is isolated. Following the proof of Proposition \ref{prop-zero-class-general}, we can define the isomorphism $\gamma:\calA(0_{12})\rightarrow \calA(U),$ by taking:

  For minimal vertices:
  \begin{equation*}
    \gamma(X_1)=Y_1,\quad \gamma(X_4)=Y_4,\quad \gamma(X_{11})=Y_{11}.
  \end{equation*}
For vertices associated to chain $\calC^{(1)}:$
  \begin{equation*}\label{eq2-ex1-prop-zero-class-gen}
  \left\{
  \begin{array}{c}
      \gamma(X_2)=-\frac{1}{2} \gamma(X_1)+Y_2=-\frac{1}{2} Y_1+Y_2, \\
       \quad\\
       \gamma(X_3)=-\gamma(X_2)+Y_3=\frac{1}{2} Y_1-Y_2+Y_3,\\
       \quad\\
      \gamma(X_5)=-\gamma(X_3)+Y_5=-\frac{1}{2} Y_1+Y_2-Y_3+Y_5, \\
      \quad\\
       \gamma(X_7)=-\gamma(X_5)+Y_7=\frac{1}{2} Y_1-Y_2+Y_3-Y_5+Y_7,\\
       \quad\\
      \gamma(X_{10})=-\frac{1}{2}\gamma(X_7)+Y_{10}=-\frac{1}{4} Y_1+\frac{1}{2}Y_2-\frac{1}{2}Y_3+\frac{1}{2}Y_5-\frac{1}{2}Y_7+Y_{10}.
    \end{array}\right.
  \end{equation*}
For vertices associated to chain $\calC^{(2)},$
\begin{equation*}\label{eq3-ex1-prop-zero-class-gen}
  \left\{
  \begin{array}{c}
      \gamma(X_6)=-\frac{1}{2}\gamma(X_3)+Y_6=-\frac{1}{4} Y_1+\frac{1}{2}Y_2-\frac{1}{2}Y_3+Y_6, \\
      \quad\\
      \gamma(X_{12})=-\frac{3}{2}\gamma(X_6)+Y_{12}=\frac{3}{8} Y_1-\frac{3}{4}Y_2+\frac{3}{4}Y_3-\frac{3}{2} Y_6+Y_{12}.
    \end{array}\right.
  \end{equation*}
And finally, for vertices associated to chain $\calC_2:$
  \begin{equation*}
    \left\{\begin{array}{c}
             \gamma(X_8)=-\frac{1}{2}\gamma(X_4)+Y_8=-\frac{1}{2}Y_4+Y_8, \\
             \quad \\
             \gamma(X_{9})=-\frac{3}{2}\gamma(X_8)+Y_{9}=
             \frac{3}{4}Y_4-\frac{3}{2}Y_8+Y_{9}.
           \end{array}\right.
  \end{equation*}
\end{example}

The following Proposition provides another evidence of the veracity of Conjecture \ref{conj-sim-vs-approx}.

\begin{proposition}\label{prop2-sim-vs-approx}
  Let $U$ be a \nSLTM. If $U\sim 0_n,$ then $U\approx 0_n.$
\end{proposition}

\begin{proof}
  Assuming that $U\sim 0_n,$ but $U\neq 0_n.$ We shall prove that $U\approx 0_n$ by induction on the number $m$ of non-zero rows of $U.$

  \begin{itemize}
    \item If $m=1,$ let $(r,c)$ be the leader of the only non-zero row of $U.$ Since $U\sim 0_n,$ Theorem \ref{theo-zero-class-n-gen} implies that either $c=1$ or $\Delta_{i,c,r}^{(2)}(U)=0$ for each $1\leq i<c.$ Both cases implies that $u_{rc}$ is the only nonzero entry of $U$ and that the ETO:  $\calQ_{(r,c)}(\quad,\frac{u_{rc}}{2})$ is well-defined over $U.$ Is not difficult to see that $0_n=\calQ_{(r,c)}(U,\frac{u_{rc}}{2}),$ and therefore $U\approx 0_n.$
    \item If $m>1,$ let $r_0$ the minimal value such that the $r_0-$th row of $U$ is different from zero. Let $(r_0,c_0)$ be the leader of the $r_0-$th row. Since $U\sim 0_n,$ Theorem \ref{theo-zero-class-n-gen} implies that either $c=1$ or $\Delta_{i,c_0,r_0}^{(2)}=0$ for each $1\leq i<c_0.$ Both cases implies that $u_{r_0c_0}$ is the only nonzero entry of $U$ and that the ETO:  $\calQ_{(r_0,c_0)}(\quad,{u_{r_0c_0}}/{2})$ is well-defined over $U.$ Let $U'=\calQ_{(r_0,c_0)}(U,{u_{r_0c_0}}/{2}),$ then most of the entries $u'_{rc}$ of $U'$ are equal to the corresponding entry $u_{rc}$ of $U,$ with the exception of those in positions of the form $(r,c_0),$ with $r\geq r_0$ such that $u_{rr_0}\neq 0.$ For those cases we have the relation:
        \begin{equation}\label{eq1-prop2-sim-vs-approx}
          u'_{rc_0}=\left\{\begin{array}{cc}
                              0 & \textrm{if}\quad r=r_0 \\
                              \quad & \quad \\
                             \frac{\Delta_{c_0,r_0,r}^{(2)}(U)}{2} & \textrm{if}\quad r>r_0
                            \end{array}\right.
        \end{equation}
        (see equation \ref{eq1-def-ETO-Q}). Therefore in the aforementioned positions, we have that $u'_{rc_0}=0.$

        If $r\geq r_0$ is such that $u_{rr_0}\neq 0$ and $(r,j)$ is the leader of the $r-$th row of $U,$ then $\Delta_{r_0,j,r}^{(2)}(U)=0$ since $U\sim 0_n.$ In particular, we have that $u_{jr_0}\neq 0$ and by definition of $r_0$ we have that $j\geq r_0.$ By equation \ref{eq1-prop2-sim-vs-approx}, we conclude that $u'_{jc_0}=0$ and then $\Delta_{c_0,j,r}^{(2)}(U')=2u'_{rc_0}+u'_{rj}u'_{jc_0}=0.$ For any other triple $(i,c,r)$ with $i<c<r$ the entries $u'_{ri},u'_{rc},u'_{ci}$ of $U'$ are the same as the corresponding entries of $U,$ then it is clear that $\Delta_{i,c,r}^{(2)}(U')=\Delta_{i,c,r}^{(2)}(U')=0$ for each $(r,c)\in \calL_{U'}$ and each $1\leq i<c.$ Theorem \ref{theo-zero-class-n-gen}, implies that $U'\sim 0_n.$

  By construction $U\approx U'$ and by equation \ref{eq1-prop2-sim-vs-approx}, $U'$ is a matrix with $m-1$ non-zero rows. By induction we conclude that $U'\approx 0_n,$ and by transitivity of relation $\approx,$ we conclude $U\approx 0_n,$ as desired.
  \end{itemize}
\end{proof}

\begin{example}\label{ex3-mega-ex-chains}
  Let $U$ be the matrix of equation \ref{eq6-ex-mega-ex-chains}. As we saw in Example \ref{ex2-mega-ex-chains}, $U\sim 0_{12}.$

  Let $\Theta$ be the matrix defined as the following product of elementary matrices:

  \begin{equation}\label{eq1-ex3-mega-ex-chains}
    \Theta=E_9\cdot E_{8}\cdot E_7\cdot E_{6}\cdot E_5\cdot E_4\cdot E_3\cdot E_2 \cdot E_1
  \end{equation}
  where
  \begin{equation}
    \begin{array}{ccc}
      E_1=Q_{(2,1)}(1/2) & E_2= Q_{(3,2)}(1)& E_3= Q_{(5,3)}(1)\\
      \quad &\quad & \quad \\
      E_4=Q_{(6,3)}({1}/{2}) & E_5=Q_{(7,5)}(1) & E_6=Q_{(8,4)}(1/2) \\
      \quad &\quad & \quad \\
      E_7=Q_{(9,8)}(3/2) & E_8= Q_{(10,7)}(1/2) & E_9=Q_{(12,6)}(3/2).
    \end{array}
  \end{equation}
  One can check that the matrix $\Theta,$ defines an isomorphism $\theta:\calA(U)\rightarrow \calA(0_{12}),$ and the factorization given in equation \ref{eq1-ex3-mega-ex-chains} yields a well-defined  process of transformation via $\calQ-$operations, that transform $U$ into $0_{12}.$ This process can be deduced from the proof of Proposition \ref{prop2-sim-vs-approx}.

  Note that for the isomorphism $\gamma:\calA(0_{12})\rightarrow \calA(U),$ defined in example \ref{ex2-mega-ex-chains}, the associated matrix $\Gamma$ can also be factored as an adequate product of elementary matrices, that yields a well-defined  process of transformation via $\calQ-$operations, that transform $0_{12}$ into $U.$ We left details to the readers.
\end{example}

\subsection{A lower bound for the number of classes}\label{ssec-lower-bound}
We would like to know, how many classes can be found in \nSLTMn, for each $n$ and how this number depends (or not) on the ground field $\F.$ Let us denote $N_n$ the number of classes in $\mathcal{TM}_n$ under $\sim.$ As we observed in Subsection \ref{ssec-small-cases}, the number $N_n,$ for $n=2,3,$ is finite and do not depend on the ground field. Although, at this point we do not know if this number is always finite or if it is independent or not on the ground field, in this subsection we obtain a lower bound for $N_n,$  that is still independent on the ground field:

\begin{theorem}\label{theo-lower-bound}
  There are at least $n-1$ classes in \nSLTMn.
\end{theorem}

\begin{proof}
We divide the proof in several steps, where basically we shall show that the matrices $B_{n,1},$ $B_{n,2},\dots,$ $B_{n,n-1}$ are representatives of pairwise different classes.

\begin{description}
\item[Step 1] First note that $0_n\sim B_{n,l}$ if and only if $l=1.$ In fact, note that the matrix $B_{n,1}$ satisfies the first condition of Theorem \ref{theo-zero-class-n-gen}, then $B_{n,1}\sim 0_n.$  On the other hand, if $l\geq 2,$ then by definition of matrix $B_{n,l}$ (see Subsection \ref{ssec-notations}), the entry in position $(n,l)$ is not zero and  $\Delta_{1,l,n}^{(2)}(B_{n,l})=1\neq 0,$ therefore the matrix $B_{n,l}$ do not satisfy the necessary condition of Theorem \ref{theo-zero-class-n-gen}, and then $B_{n,l}\nsim 0_n$ as desired.
\item[Step 2]  Let $T=B_{n,l}=[t_{ij}]$ and $S=B_{n,m}=[s_{ij}]$ with $1< l<m<n.$ If we assume that there is an isomorphism $\gamma:\calA(T)\rightarrow \calA(S),$ then by Theorem \ref{coro-Key-condition}, we have
  \begin{equation}\label{Key-EQ-adapted3}
    2\gamma_{ir}\gamma_{kr}+\gamma_{kr}^2s_{ki}=
    \sum_{j<r}t_{rj}(\gamma_{kj}\gamma_{kr}s_{ki}+\gamma_{kj}\gamma_{ir}+\gamma_{ij}\gamma_{kr}),
    \quad(1\leq r\leq n, 1\leq i<k\leq n).
  \end{equation}
  \begin{description}
 \item[Step 2.1]  In particular if we fix a $r<n,$ then for each pair $1\leq i<k<n,$ equation \ref{Key-EQ-adapted3} reduces to
  \begin{equation}\label{Key-EQ-adapted4}
    2\gamma_{ir}\gamma_{kr}=0.
  \end{equation}
 By our assumption, necessarily for each $r$ there is a $1\leq c_r\leq n$ such that $\gamma_{c_rr}\neq 0,$ (Since the associated matrix $\Gamma$ has to be invertible). If we assume that $c_r<n,$ then equation \ref{Key-EQ-adapted4} implies that $\gamma_{ir}=0,$ whenever $1\leq i<n$ and $i\neq c_r.$  On the other hand, if we assume that $c_r=n,$ and we replace $k=n$ in equation \ref{Key-EQ-adapted3}, we obtain $\gamma_{nr}s_{ni}=0,$ and therefore $s_{ni}=0,$ for all $1\leq i<n,$ but this is a contradiction with the definition of the matrix $S=B_{n,m}.$
 The last analysis implies that for each $1\leq r<n,$ there is a $1\leq c_r<n$ such that $\gamma_{c_rr}\neq 0$ and $\gamma_{ir}=0$ for each $1\leq i<n,$ with $i\neq c_r.$

 Now we replace in equation \ref{Key-EQ-adapted3}, $k=n$ and $1\leq i\leq m$ such that $i\neq c_r,$ then we obtain $\gamma_{nr}^2s_{ni}=0.$ Since $s_{ni}=1$ by definition of the matrix $S=B_{n,m},$ then we conclude that $\gamma_{nr}=0.$ That is valid for each $1\leq r<n.$ Therefore, since $\det(\Gamma)\neq 0,$ we conclude that $\gamma_{nn}\neq 0.$

 By the last analysis we also conclude that for each pair $1\leq r_1<r_2<n,$ we have $c_{r_1}\neq c_{r_2}.$ In particular, we conclude that there are exactly $m$ indices $r\in\{1,\dots,n-1\}$ such that $c_r\leq m.$

\item[Step 2.2] Let us fix now a pair $1\leq r_1\neq r_2<n$ such that $c_{r_1}<c_{r_2}.$ If we replace $r=n,i=c_{r_1},k=c_{r_2}$ in equation \ref{Key-EQ-adapted3}, we obtain

  \begin{equation}\label{Key-EQ-adapted50}
    2\gamma_{c_{r_1}n}\gamma_{c_{r_2}n}+\gamma_{c_{r_2}n}^2s_{c_{r_2}c_{r_1}}=
    \sum_{j < n}t_{nj}(\gamma_{c_{r_2}j}\gamma_{c_{r_2}n}s_{c_{r_2}c_{r_1}}+\gamma_{c_{r_2}j}\gamma_{c_{r_1}n}
    +\gamma_{c_{r_1}j}\gamma_{c_{r_2}n}).
  \end{equation}

  and by definition of the matrices $T=B_{n,l},S=B_{n,m}$  and the indices $c_{r_1},c_{r_2},$ we have that $s_{c_{r_2}c_{r_1}}=0,$ $\gamma_{c_{r_1}j}=0$ whenever $j\neq r_1$ and  $\gamma_{c_{r_2}j}=0$ whenever $j\neq r_2$ , therefore equation \ref{Key-EQ-adapted50} reduces to

  \begin{equation}\label{Key-EQ-adapted5}
    2\gamma_{c_{r_1}n}\gamma_{c_{r_2}n}=
    t_{nr_2}\gamma_{c_{r_2}r_2}\gamma_{c_{r_1}n}
    +t_{nr_1}\gamma_{c_{r_1}r_1}\gamma_{c_{r_2}n}.
  \end{equation}

  Note that by the symmetry of equation \ref{Key-EQ-adapted5}, we can omit the hypothesis that $c_{r_1}<c_{r_2}.$

\item[Step 2.3] If we replace $k=r=n$ in equation  \ref{Key-EQ-adapted3}, we obtain
 \begin{equation}\label{Key-EQ-adapted51}
    2\gamma_{in}\gamma_{nn}+\gamma_{nn}^2s_{ni}=
    \sum_{j<n}t_{nj}(\gamma_{nj}\gamma_{nr}s_{ni}+\gamma_{nj}\gamma_{in}+\gamma_{ij}\gamma_{nn}),
    \quad(1\leq i< n).
  \end{equation}
  By definition of the matrix $T=B_{n,l}$ (in particular using that $t_{nj}=1$ if $j\leq l$ and $t_{nj}=0$ if $l<j<n$), and by the facts that $\gamma_{nj}=0$ for each $1\leq j<n,$ and $\gamma_{nn}\neq 0,$ we obtain, after simplification by $\gamma_{nn},$ that
  \begin{equation}\label{Key-EQ-adapted52}
    2\gamma_{in}+\gamma_{nn}s_{ni}=
    \sum_{j\leq l}\gamma_{ij},\quad(1\leq i< n).
  \end{equation}
\item[Step 2.4] Let us fix $1\leq r <n.$ If we replace $i=c_r$ in equation \ref{Key-EQ-adapted52},  then we have that $\sum_{j\leq l}\gamma_{ij}=\gamma_{c_{r}r},$ and then we obtain

  \begin{equation}\label{Key-EQ-adapted6}
    \gamma_{c_rn}=\left\{\begin{array}{cc}
    \frac{1}{2}(\gamma_{c_rr}-\gamma_{nn}s_{nc_r}) & \quad\textrm{if}\quad r\leq l \\
    \quad&\quad\\
    -\frac{1}{2}\gamma_{nn}s_{nc_r} & \quad\textrm{otherwise}
    \end{array}\right.
  \end{equation}
\item[Step 2.5]  Let us fix $1\leq r_1< r_2<n$ with $r_2>l,$ such that $c_{r_1},c_{r_2}\leq m$ (we know there is at least two of them satisfying that condition), in particular $s_{nc_{r_1}}=s_{nc_{r_2}}=1.$ Therefore we have the following two cases:

  \begin{enumerate}
  \item If $r_1\leq l<r_2$ then equation \ref{Key-EQ-adapted5} turns into
    \begin{equation}\label{Key-EQ-adapted61}
    2\gamma_{c_{r_1}n}\gamma_{c_{r_2}n}=
    \gamma_{c_{r_1}r_1}\gamma_{c_{r_2}n}
  \end{equation}
  since $t_{nr_1}=1$ and $t_{nr_2}=0.$ Then by equations \ref{Key-EQ-adapted6} and \ref{Key-EQ-adapted61}, we obtain
  \begin{equation*}
   -\frac{1}{2}(\gamma_{c_{r_1}r_1}-\gamma_{nn})
    \gamma_{nn}=
    -\frac{1}{2}\gamma_{c_{r_1}r_1}\gamma_{nn},
  \end{equation*}
   from where we obtain $\gamma_{nn}=0,$ which is a contradiction.
  \item If $l<r_1,r_2,$ then equation \ref{Key-EQ-adapted5} turns into
    \begin{equation}\label{Key-EQ-adapted62}
    2\gamma_{c_{r_1}n}\gamma_{c_{r_2}n}=
    0.
  \end{equation}
  since in this case $t_{nr_1}=t_{nr_2}=0.$ But equation \ref{Key-EQ-adapted6} and our assumption that $c_{r_1},c_{r_2}\leq m$ imply that:
  \begin{equation*}
    \gamma_{c_{r_1}n}=\gamma_{c_{r_2}n}=-\frac{1}{2}\gamma_{nn}\neq 0,
  \end{equation*}
  which is a contradiction.
  \end{enumerate}
The last two cases imply that whenever $c_{r}\leq m$ we necessarily have $r\leq l.$ Then there are at most $l$ indices $1\leq r<n$ such that $c_{r}\leq m,$ but this is a contradiction, since we know from \textbf{Step 2.1}, that there are exactly $m$ of them satisfying that condition (recall that $l<m$).
\end{description}
By the last analysis we conclude that $\calA(T)$ and $\calA(S)$ are not isomorphic, and then $B_{n,l}\nsim B_{n,m},$ whenever $1<l<m<n.$
\item[Step 3] By \textbf{Step 1} and \textbf{Step 2}, we conclude that the matrices $B_{n,1},B_{n,2},\dots,B_{n,n-1}$ are representatives of pairwise different classes under relation $\sim.$ Therefore there are at least $n-1$ different classes in \nSLTMn.
\end{description}
\end{proof}

Motivated from the results of subsection \ref{ssec-small-cases} and \ref{ssec-lower-bound}, we open the following questions:

\begin{description}
  \item[Q1:]  Is the number $N_n$ always finite?
 \item[Q2:] Is the number $N_n$ independent of the ground field?
\end{description}
We do not have an answer for this questions yet, but motivated by our intuition, we define the following:

\begin{conjecture}\label{conj-on-Nn}
  The number $N_n$ for $n>3$ depends on the ground field. Moreover the number $N_n$ is finite if and only if the ground field is finite.
\end{conjecture}

We expect to be able to prove or disprove Conjecture \ref{conj-on-Nn} in future works.


\subsection{$JM-$subalgebras in Soergel Calculus.}
In this subsection we consider the field of complex numbers $\F=\mathbb{C}.$
Given a Coxeter system of finite rank $(W,S),$ (with $S=\{s_a:a\in J\}$ for certain set of indices $J$), we consider the \emph{geometric realization} $(\mathfrak{h},\{\alpha^{\vee}_a\},\{\alpha_a\})$ of $(W,S)$ over the field $\mathbb{C}.$ That is a $\mathbb{C}-$vector space $\mathfrak{h},$ and subsets $\{\alpha^{\vee}_a\}\subset\mathfrak{h},$ and  $\{\alpha_a\}\subset\mathfrak{h}^{\ast}=\textrm{Hom}_{\mathbb{C}}(\mathfrak{h},\mathbb{C}),$
where
\begin{equation}\label{definition-action-root-coroot}
  \langle \alpha_a,\alpha^{\vee}_b \rangle:=\left\{ \begin{array}{cc}
                                                     -2\cos(\pi/{\mathbf{m}(a,b)}) & \quad\textrm{if}\quad \mathbf{m}(a,b)<\infty \\
                                                     -2 & \quad\textrm{if}\quad \mathbf{m}(a,b)=\infty
                                                   \end{array}\right.
\end{equation}

and the equations:


\begin{equation}\label{definition-geometric-representation}
  \left\{\begin{array}{cc}
           \beta\cdot s_{a} =\beta-\langle\alpha_a,\beta \rangle\alpha^{\vee}_a, & \quad (\beta \in \mathfrak{h}). \\
           \quad & \quad \\
           \gamma\cdot s_{a}=\gamma-\langle\gamma,\alpha^{\vee}_a\rangle \alpha_a, & \quad (\gamma \in \mathfrak{h}^{\ast}).
         \end{array}\right.
\end{equation}
define \emph{representations} of the group $W$ on $\mathfrak{h}$ and $\mathfrak{h}^{\ast}$ respectively (see \cite{H1} and \cite{EW} for details).

For any $\gamma\in \textrm{Span}(\{\alpha_a:a\in J\})\subset\mathfrak{h}^{\ast}$ and any $a\in J$ we denote $\partial_a(\gamma)=\langle \gamma,\alpha^{\vee}_a\rangle.$
Let $W^{\ast}$ the free monoid generated by the symbols $\{s_a:a\in J\}.$  In \cite{Lobos2}, for each element $\ulu=s_{a_1}\cdots s_{a_p}\in W^{\ast}$ we define the matrix:
\begin{equation*}
\mathbf{T}_{\ulu}=\left[\begin{array}{ccccc}
                                   0 & 0 & 0 & \cdots& 0 \\
                                   \partial_{a_1}(\alpha_{a_2}) & 0 & 0 & \cdots& 0 \\
                                   \partial_{a_1}(\alpha_{a_3}\cdot s_{a_2}) & \partial_{a_2}(\alpha_{a_3}) & 0 & \cdots& 0 \\
                                   \vdots & \vdots & \ddots & \ddots& \vdots\\
                                   \partial_{a_1}(\alpha_{a_p}\cdot s_{a_{p-1}}\dots s_{a_2}) & \partial_{a_2}(\alpha_{a_p}\cdot s_{a_{p-1}}\cdots s_{a_{3}}) & \cdots& \partial_{a_{p-1}}(\alpha_{a_p}) & 0 \\
                                 \end{array}\right]
\end{equation*}
  and the nil graded algebra $\calA_{\ulu}=\calA(\mathbf{T}_{\ulu})$ associated to $\mathbf{T}_{\ulu}.$ We have proved in \cite{Lobos2}, that the algebra $\calA_{\ulu}$ corresponds to the Gelfand-Tsetlin subalgebra $\bcalJ_{\ulu}$ of the endomorphism algebra $\textrm{End}(\ulu)$ in the diagrammatic Soergel Category of Elias and Williamson (see \cite{EW}), relative to the Libedynski's light leaves basis (see \cite{LibLightLeaves}) and the set of Jucys-Murphy elements defined by Ryom-Hansen in \cite{Steen2}.

  In \cite{Lobos2} we show that for the case of a Coxeter system of type $\widetilde{A},$ the algebras $\calA_{\ulu}$ are invariant under commuting moves, and we proved that by defining explicit isomorphisms. Now in theorem \ref{theo-appl-commuting-moves} we generalize this fact for any other type. In this case we use the techniques developed in this article, to obtain a cleaner proof.

  Let us recall first, certain technical definitions coming from \cite{Lobos2}:
  
  Given a \nSLTM, $T,$  a \mSLTM, $S$ and a $m\times n-$matrix $C,$ we define the matrix $T\nabla_{C}S$ by blocks, as follows:
    \begin{equation*}
      T\nabla_{C}S=\left[\begin{matrix}
                           T & 0 \\
                           C & S
                         \end{matrix}\right]
    \end{equation*}

  \begin{definition}
  \begin{enumerate}
    \item Let $\ulv=s_{b_1}\cdots s_{b_{l}}\in W^{\ast}$ we define the \emph{column vector} associated to $\ulv$ as the matrix $Q_{\ulv}$ given by:
        \begin{equation*}
          Q_{\ulv}=\left[\begin{matrix}
                           q_1{(\ulv)} \\
                           q_2{(\ulv)} \\
                           \vdots\\
                           q_l{(\ulv)}
                         \end{matrix}\right]=\left[\begin{matrix}
                           \alpha_{b_1} \\
                           \alpha_{b_2}\cdot s_{b_1} \\
                           \vdots\\
                           \alpha_{b_{l}}\cdot (s_{b_{l-1}}\cdots s_{b_1})\end{matrix}\right],
        \end{equation*}
    \item Given two elements of $W^{\ast},$ say $\ulu=s_{a_1}\cdots s_{a_{p}}$ and $\ulv=s_{b_1}\cdots s_{b_{l}}$ we define the matrix $\bpartial_{\ulu}(Q_{\ulv})$ by
    \begin{equation*}
    \bpartial_{\ulu}(Q_{\ulv})= \left[\begin{array}{cccc}
                              \partial_{a_1}(q_{1}(\ulv)\cdot \ulu_{>1}), & \partial_{a_2}(q_{1}(\ulv)\cdot \ulu_{>2}), & \cdots
                              &,\partial_{a_p}(q_{1}(\ulv)) \\
                              \partial_{a_1}(q_{2}(\ulv)\cdot \ulu_{>1}), & \partial_{a_2}(q_{2}(\ulv)\cdot \ulu_{>2}), &\cdots &,\partial_{a_p}(q_{2}(\ulv))\\
                              \vdots & \vdots & \vdots & \vdots\\
                             \partial_{a_1}(q_{l}(\ulv)\cdot \ulu_{>1}), & \partial_{a_2}(q_{l}(\ulv)\cdot \ulu_{>2}), & \cdots &,\partial_{a_p}(q_{l}(\ulv)) \\
                          \end{array}\right].
    \end{equation*}
    where $\ulu_{>r}=s_{a_{r+1}}\cdots s_{p}.$
  \end{enumerate}
  \end{definition}

  In \cite{Lobos2}, we show that for any composed expression $\ulw=\ulu\ulv,$ we have that:
  \begin{equation*}
    \mathbf{T}_{\ulw}=\mathbf{T}_{\ulu}\nabla_{C}\mathbf{T}_{\ulv},\quad
    \textrm{where}\quad C=\bpartial_{\ulu}(Q_{\ulv}).
  \end{equation*}
  Moreover in \cite{Lobos2}, we proved that for $\ulz=\ulu\ulv\ulw$ we have the following (kind of) associative property:
  \begin{equation*}
    \mathbf{T}_{\ulz}=\left(\mathbf{T}_{\ulu}\nabla_{C}\mathbf{T}_{\ulv}\right)\nabla_{D}
    \mathbf{T}_{\ulw}=\mathbf{T}_{\ulu}\nabla_{C'}\left(\mathbf{T}_{\ulv}\nabla_{D'}\mathbf{T}{\ulw}\right)
  \end{equation*}
  where
  \begin{equation*}
    C=\bpartial_{\ulu}(Q_{\ulv}), D=\bpartial_{\ulu\ulv}(Q_{\ulw}),
    \quad C'=\bpartial_{\ulu}(Q_{\ulv\ulw}), D'=\bpartial_{\ulv}(Q_{\ulw}).
  \end{equation*}
Now we are ready to prove the following:
\begin{theorem}\label{theo-appl-commuting-moves}
  For any Coxeter system $(W,S),$ the algebra $\calA_{\ulu}$ is invariant under commuting moves.
\end{theorem}

\begin{proof}
  Let $b,c\in J$ such that $\mathbf{m}(b,c)=2,$ that its, the elements $s_b$ and $s_c$ of the group $W$ commute. Let $\ulu=s_{a_1}\cdots s_{a_{p-1}}s_bs_cs_{a_{p+2}}\cdots s_{a_{l}}$ and $\ulv=s_{a_1}\cdots s_{a_{p-1}}s_cs_bs_{a_{p+2}}\cdots s_{a_{l}},$ then $\ulu,\ulv$ are elements of $W^{\ast}$ that differ by a commuting move.
  \begin{enumerate}
    \item Let $\ulw=s_{a_1}\cdots s_{a_{p-1}}$ and $\ulz=s_{a_{p+2}}\cdots s_{a_{l}}.$ Then $\ulu=\ulw\cdot s_{b}s_{c}\cdot \ulz$ and $\ulv=\ulw \cdot s_cs_b \cdot\ulz.$  Therefore:
        \begin{equation}\label{eq1-theo-appl-commuting-moves}
          \mathbf{T}_{\ulu}=\mathbf{T}_{\ulw}\nabla_X\left(\mathbf{T}_{s_bs_c}\nabla_{Y} \mathbf{T}_{\ulz}\right),\quad
          \mathbf{T}_{\ulv}=\mathbf{T}_{\ulw}\nabla_{X'}\left(\mathbf{T}_{s_cs_b}\nabla_{Y'} \mathbf{T}_{\ulz}\right)
        \end{equation}
        where
        \begin{equation*}
         X=\bpartial_{\ulw}(Q_{s_bs_c\ulz}),Y=\bpartial_{s_bs_c}(Q_{\ulz}),
         \quad X'=\bpartial_{\ulw}(Q_{s_cs_b\ulz}),Y'=\bpartial_{s_cs_b}(Q_{\ulz}).
        \end{equation*}

    \item Note that by definition $\partial_{b}(\alpha_c)=0=\partial_{c}(\alpha_b),$ $\alpha_c\cdot s_b=\alpha_c$ and $\alpha_b\cdot s_c=\alpha_b.$ Therefore we have:
        \begin{equation*}
          \mathbf{T}_{s_bs_c}=\left[\begin{matrix}
                                      0 & 0 \\
                                      0 & 0
                                    \end{matrix}\right]=\mathbf{T}_{s_cs_b},\quad
        Q_{s_bs_c}=\left[\begin{matrix}
                           \alpha_b \\
                           \alpha_c
                         \end{matrix}\right],\quad Q_{s_cs_b}=\left[\begin{matrix}
                           \alpha_c \\
                           \alpha_b
                         \end{matrix}\right].
        \end{equation*}
        By the same reason, we have that two first rows of $X$ are given by:
        \begin{equation*}
          X_1=\left[\begin{matrix}
                  \partial_{a_1}(\alpha_b\cdot s_{a_{p-1}}\cdots s_{a_2}) & \partial_{a_2}(\alpha_b\cdot s_{a_{p-1}}\cdots s_{a_3}) & \cdots & \partial_{a_{p-1}}(\alpha_b) \\
                  \partial_{a_1}(\alpha_c\cdot s_{a_{p-1}}\cdots s_{a_2}) & \partial_{a_2}(\alpha_c\cdot s_{a_{p-1}}\cdots s_{a_3}) & \cdots & \partial_{a_{p-1}}(\alpha_c)
                \end{matrix}\right]
        \end{equation*}
       analogously the two first rows of $X'$ are given by:
        \begin{equation*}
          X_1'=\left[\begin{matrix}
                 \partial_{a_1}(\alpha_c\cdot s_{a_{p-1}}\cdots s_{a_2}) & \partial_{a_2}(\alpha_c\cdot s_{a_{p-1}}\cdots s_{a_3}) & \cdots & \partial_{a_{p-1}}(\alpha_c) \\
                  \partial_{a_1}(\alpha_b\cdot s_{a_{p-1}}\cdots s_{a_2}) & \partial_{a_2}(\alpha_b\cdot s_{a_{p-1}}\cdots s_{a_3}) & \cdots & \partial_{a_{p-1}}(\alpha_b)
                \end{matrix}\right]
        \end{equation*}

    \item Note that for every $\gamma\in \textrm{Span}(\{\alpha_a:a\in J\}),$ we have that $\partial_b(\gamma\cdot s_c)=\partial_b(\gamma).$ Let us denote
          \begin{equation*}
            \gamma_{j}=\left\{\begin{array}{cc}
                                \alpha_{a_j}\cdot s_{a_{j-1}}\cdots s_{a_{p+2}}, & \textrm{if}\quad p+2<j\leq l \\
                                \alpha_{a_{p+2}} & \textrm{if}\quad j=p+2
                              \end{array}\right.
          \end{equation*}
          then
            \begin{equation*}
              Y=\left[\begin{matrix}
                      \partial_{b}(\gamma_{p+2}) & \partial_{c}(\gamma_{p+2}) \\
                      \vdots & \vdots \\
                      \partial_{b}(\gamma_{l}) & \partial_{c}(\gamma_{l})
                    \end{matrix}\right],\quad
              Y'=     \left[\begin{matrix}
                      \partial_{c}(\gamma_{p+2}) & \partial_{b}(\gamma_{p+2}) \\
                      \vdots & \vdots \\
                      \partial_{c}(\gamma_{l}) & \partial_{b}(\gamma_{l})
                    \end{matrix}\right]
            \end{equation*}

           \item Following with the notation of the $\gamma'$s, we have that the rows of the matrices $X,X'$ from the third to the last one are given respectively by:
            \begin{equation*}
              X_2=\left[\begin{matrix}
                      \partial_{a_1}((\gamma_{p+2}\cdot s_bs_c)\cdot s_{a_p}\cdots s_{a_2}) & \cdots & \partial_{a_{p-1}}(\gamma_{p+2}\cdot s_bs_c)  \\
                      \partial_{a_1}((\gamma_{p+3}\cdot s_bs_c)\cdot s_{a_p}\cdots s_{a_2}) & \cdots & \partial_{a_{p-1}}(\gamma_{p+3}\cdot s_bs_c)  \\
                      \vdots&\cdots&\vdots\\
                      \partial_{a_1}((\gamma_{l}\cdot s_bs_c)\cdot s_{a_p}\cdots s_{a_2}) & \cdots & \partial_{a_{p-1}}(\gamma_{l}\cdot s_bs_c)
                    \end{matrix}\right]
            \end{equation*}
           and
           \begin{equation*}
              X_2'=\left[\begin{matrix}
                      \partial_{a_1}((\gamma_{p+2}\cdot s_cs_b)\cdot s_{a_p}\cdots s_{a_2}) & \cdots & \partial_{a_{p-1}}(\gamma_{p+2}\cdot s_bs_c)  \\
                      \partial_{a_1}((\gamma_{p+3}\cdot s_cs_b)\cdot s_{a_p}\cdots s_{a_2}) & \cdots & \partial_{a_{p-1}}(\gamma_{p+3}\cdot s_bs_c)  \\
                      \vdots&\cdots&\vdots\\
                      \partial_{a_1}((\gamma_{l}\cdot s_cs_b)\cdot s_{a_p}\cdots s_{a_2}) & \cdots & \partial_{a_{p-1}}(\gamma_{l}\cdot s_bs_c)
                    \end{matrix}\right]
            \end{equation*}
            But since $s_b$ and $s_c$ commute, we clearly have that $X_2=X_2'.$
           \item By equation \ref{eq1-theo-appl-commuting-moves} we have that:
           \begin{equation*}
             \mathbf{T}_{\ulu}=\left[\begin{matrix}
                                       \mathbf{T}_{\ulw} & 0 & 0 \\
                                       X_1 &\mathbf{T}_{s_bs_c} & 0 \\
                                       X_2 & Y & \mathbf{T}_{\ulz}
                                     \end{matrix}\right],\quad
             \mathbf{T}_{\ulv}=\left[\begin{matrix}
                                       \mathbf{T}_{\ulw} & 0 & 0 \\
                                       X_1' &\mathbf{T}_{s_cs_b} & 0 \\
                                       X_2' & Y' & \mathbf{T}_{\ulz}
                                     \end{matrix}\right],
           \end{equation*}
    and by the analysis above, we conclude that $\mathbf{T}_{\ulu}\sim \mathbf{T}_{\ulv}.$ Moreover $\mathbf{T}_{\ulv}=\calF_{(p,p+1)}(\mathbf{T}_{\ulu}),$ and therefore (by lemma \ref{Lemma-ETO-F} ) the algebras $\calA_{\ulu}$ and $\calA_{\ulv}$ are isomorphic.
  \end{enumerate}
\end{proof}

Theorem \ref{theo-appl-commuting-moves} naturally opens the following question:

\begin{description}
  \item[Q3:]  Is the algebra $\calA_{\ulu}$ invariant under other braid moves?
\end{description}
Fortunately we have an answer for \textbf{Q3}, although it is not a positive one:
\begin{proposition}\label{prop-calA-not-invariant-under-braid-moves}
  The algebra $\calA_{\ulu}$ is not invariant under braid moves.
\end{proposition}
\begin{proof}
  It is enough to provide an example where two equivalent by braid moves expressions $\ulu$ and $\ulv$ generate two not isomorphic algebras $\calA_{\ulu}$ and $\calA_{\ulv}.$

  Let $(W,S)$ be the Coxeter system of type $\widetilde{A}_3.$ Let $\ulu=s_0s_1s_2s_3s_2$ and $\ulv=s_0s_1s_3s_2s_3$ then:
  \begin{enumerate}
    \item It is clear that $\ulu$ and $\ulv$ are equivalent by braid moves.
    \item The corresponding matrices $\mathbf{T}_{\ulu}$ and $\mathbf{T}_{\ulv}$ are given by:
        \begin{equation*}
          \mathbf{T}_{\ulu}=\left[\begin{matrix}
                              0 & 0 & 0 & 0 & 0 \\
                              -1 & 0 & 0 &0 & 0 \\
                              -1 & {\color{red}0} & 0 & 0 & 0 \\
                              -2 & -1 & -1 & 0 & 0 \\
                              -1 & {\color{blue}-1} & 1 & -1 & 0
                            \end{matrix}\right],\quad
          \mathbf{T}_{\ulv}=\left[\begin{matrix}
                              0 & 0 & 0 & 0 & 0 \\
                              -1 & 0 & 0 &0 & 0 \\
                              -1 & {\color{blue}-1} & 0 & 0 & 0 \\
                              -2 & -1 & -1 & 0 & 0 \\
                              -1 & {\color{red}0} & 1 & -1 & 0
                            \end{matrix}\right].
        \end{equation*}
        \item Considering the matrices $T_1,T_2,T$ defined by:
        \begin{equation*}
          T_1=\calQ_{(2,1)}\left(\mathbf{T}_{\ulu},{-\frac{1}{2}}\right),\quad T_2=\calQ_{(3,1)}\left(T_1,{-\frac{1}{2}}\right),\quad
          T=\calP_{4}\left(T_2,-1\right)
        \end{equation*}
        one can check that they are well defined (see subsection for details) and by definition $\mathbf{T}_{\ulu}\approx T_1\approx T_2\approx T.$ More specifically we have:
        \begin{equation*}
          \mathbf{T}_{\ulu} \approx T=
          \left[\begin{matrix}
          0 & 0 & 0 & 0 & 0 \\
          {0} & 0 & 0 &0 & 0 \\
          {0} & 0 & 0 & 0 & 0 \\
          {1} & {1} & {1} & 0 & 0 \\
          {-1} & -1 & 1 & {1} & 0
          \end{matrix}\right].
        \end{equation*}

        \item  Let $t_{ij}$ be the entries of $T,$ let $s_{ij}$ be the entries of the matrix $S=\mathbf{T}_{\ulv}.$
         By theorem \ref{theo-sim-vs-approx}, if the algebras $\calA_{\ulu}$ and $\calA_{\ulv}$ were isomorphic, then necessarily there will be an isomorphism $\gamma:\calA(T)\rightarrow\calA(S).$ Assuming that, then the associated matrix $\Gamma=[\gamma_{ij}]$ would satisfy equation \ref{Key-EQ-coro} of corollary \ref{coro-Key-condition}. In particular we would have for $1\leq r\leq 3$ and any pair $1\leq i<k\leq 5$ that:
        \begin{equation}\label{eq1-prop-calA-not-invariant-under-braid-moves}
          (2\gamma_{ir}+\gamma_{kr}s_{ki})\gamma_{kr}=0.
        \end{equation}
        \item Equation \ref{eq1-prop-calA-not-invariant-under-braid-moves} in particular implies that:

  \begin{equation}\label{eq2-prop-calA-not-invariant-under-braid-moves}
    \left\{\begin{array}{cc}
             (2\gamma_{4r}-\gamma_{5r})\gamma_{5r}=0, & \quad (k=5,i=4) \\
             \quad & \quad \\
             (2\gamma_{3r}+\gamma_{5r})\gamma_{5r}=0, & \quad(k=5,i=3) \\
             \quad & \quad \\
            (2\gamma_{3r}-\gamma_{4r})\gamma_{4r}=0, & \quad(k=4,i=3)
           \end{array}\right.
  \end{equation}
  If we assume that $\gamma_{5r}\neq 0,$ then system \ref{eq2-prop-calA-not-invariant-under-braid-moves} implies that $\gamma_{4r}\neq 0$ and $\gamma_{3r}\neq 0.$ Moreover we obtain:
  \begin{equation}\label{eq3-prop-calA-not-invariant-under-braid-moves}
    \left\{\begin{array}{cc}
             \gamma_{5r}=2\gamma_{4r} \\
             \gamma_{5r}=-2\gamma_{3r}\\
            \gamma_{4r}=2\gamma_{3r}
           \end{array}\right.
  \end{equation}
  but the system \ref{eq3-prop-calA-not-invariant-under-braid-moves} implies that $3\gamma_{4r}=0,$ and this is a contradiction, since the characteristic of $\F$ is different from $3$ and $\gamma_{4r}\neq0.$ Therefore we conclude that $\gamma_{5r}=0.$
  \item Equation \ref{eq1-prop-calA-not-invariant-under-braid-moves} implies that:
  \begin{equation}\label{eq4-prop-calA-not-invariant-under-braid-moves}
     \left\{\begin{array}{cc}
             (2\gamma_{3r}-\gamma_{4r})\gamma_{4r}=0, & \quad (k=4,i=3) \\
             \quad & \quad \\
             (2\gamma_{2r}-\gamma_{4r})\gamma_{4r}=0, & \quad(k=4,i=2) \\
             \quad & \quad \\
            (2\gamma_{2r}-\gamma_{3r})\gamma_{3r}=0, & \quad(k=3,i=2)
           \end{array}\right.
  \end{equation}
  If we assume that $\gamma_{4r}\neq 0,$ then system \ref{eq4-prop-calA-not-invariant-under-braid-moves} implies that $\gamma_{3r}\neq0$ and $\gamma_{2r}\neq 0.$ Moreover we obtain:
  \begin{equation}\label{eq5-prop-calA-not-invariant-under-braid-moves}
    \left\{\begin{array}{cc}
             \gamma_{4r}=2\gamma_{3r} \\
             \gamma_{4r}=2\gamma_{2r}\\
            \gamma_{3r}=2\gamma_{2r}
           \end{array}\right.
  \end{equation}
  but system \ref{eq5-prop-calA-not-invariant-under-braid-moves} implies that $2\gamma_{2r}=0,$ but this is a contradiction since the characteristic of $\F$ is different from $2$ and $\gamma_{2r}\neq 0.$ Therefore we conclude that $\gamma_{4r}=0.$
  \item Let $\{X_1,X_2,X_3,X_4,X_5\}$ and $\{Y_1,Y_2,Y_3,Y_4,Y_5\}$ be the sets of generators of the algebras $\calA(T)$ and $\calA(S)$ respectively. From the last two items, we conclude that if there exists an isomorphism $\gamma:\calA(T)\rightarrow\calA(S),$ then necessarily $\gamma(X_i)\in\textrm{Span}(Y_1,Y_2,Y_3)$ for $i=1,2,3.$ In particular this implies that the matrices $T|_3$ and $S|_3$ are in the same class under relation $\sim,$ but this contradicts theorem \ref{theo-zero-class-n-gen}. In fact
      \begin{equation*}
        T|_3=0_3=\left[\begin{matrix}
                     0 & 0 & 0 \\
                     0 & 0 & 0 \\
                     0 & 0 & 0
                   \end{matrix}\right],\quad S|_3=\left[\begin{matrix}
                     0 & 0 & 0 \\
                     -1 & 0 & 0 \\
                     -1 & -1 & 0
                   \end{matrix}\right],\quad\textrm{and}\quad \Delta_{1,2,3}^{(2)}(S)=-1\neq 0.
      \end{equation*}
  \end{enumerate}
  The analysis above implies that $\calA_{\ulu}$ and $\calA_{\ulv}$ are not isomorphic.
\end{proof}

The last Proposition opens the door for a new question:

\begin{description}
  \item[Q4:] Is there any other \emph{move} on the expression $\ulu$ that left the algebra $\calA_{\ulu}$ invariant.
\end{description}
We leave \textbf{Q4} open. For now, it is clear that new developments and techniques in the $\bcalNT$ category are needed to obtain a complete classification of the numerous $JM-$subalgebras in Soergel calculus. We hope to address this topic in future work.
\medskip
\section{Declarations}
\textbf{Acknowledgement}
\medskip

It is a pleasure to thank Steen Ryom-Hansen for stimulating conversations on Jucys-Murphy elements and $JM-$subalgebras.

\medskip
\textbf{Funding:}
This work was supported by:
\begin{itemize}
  \item {ANID, FONDECYT de Postdoctorado 2024}, N°3240046.
  \item {ANID, Subvención a la Instalación en la Academia 2024}, N°85240053.
\end{itemize}

\medskip

\textbf{Data Availability Statement:}

Data sharing not applicable to this article as no datasets were generated or analysed during the current study.
\medskip

\textbf{Competing interest:}
The author have no competing interest to declare that are relevant to the content of this article.
\medskip

\textbf{Author contribution:}
The author of this article is the only person that develop, write and review this manuscript.




\sc
diego.lobosm@uv.cl, Universidad de Valpara\'iso, Chile.

\end{document}